\documentclass[twocolumn]{arxiv}

\usepackage{soul} 
\usepackage{caption}

\usepackage{amssymb}
\usepackage{amsmath}
\usepackage{graphicx}
\graphicspath{ {fig/} }
\usepackage{epstopdf}
\usepackage{hyperref}

\begin{document}

\begin{frontmatter}

\title{Distributed methods for synchronization of orthogonal matrices over graphs\thanksref{footnoteinfo}}

\thanks[footnoteinfo]{The authors gratefully acknowledge the financial support from the Fonds
National de la Recherche, Luxembourg (FNR8864515, FNR6538106).}

\author[Johan]{Johan Thunberg}
\ead{johan.thunberg@uni.lu},
\author[Johan,Florian1]{Florian Bernard}\ead{bernard.florian@uni.lu},
\author[Johan,Cambrdige]{Jorge Goncalves}
\ead{jorge.goncalves@uni.lu}

\address[Johan]{Luxembourg Centre for Systems Biomedicine, University of Luxembourg, Esch-sur-Alzette, LUXEMBOURG}

\address[Florian1]{
Centre Hospitalier de Luxembourg, Luxembourg City, LUXEMBOURG}

\address[Cambrdige]{Control Group, Department of Engineering, University of Cambridge, Cambridge, UNITED KINGDOM}


\begin{abstract}
This paper addresses the problem of
synchronizing orthogonal matrices over directed graphs. For synchronized transformations (or matrices), composite transformations over loops equal the identity. We formulate the synchronization problem as a least-squares optimization problem with nonlinear constraints. The synchronization problem appears as one of the key components in applications ranging from 3D-localization to image registration. The main contributions of this work can be summarized as the introduction of two novel algorithms; one for symmetric graphs and one for graphs that are possibly asymmetric.
Under general conditions, the former has guaranteed convergence to the solution of a spectral relaxation to the synchronization problem. The latter is stable
for small step sizes when the graph is quasi-strongly connected.
The proposed methods 
are verified in numerical simulations.
\end{abstract}

\end{frontmatter}

\section{Introduction}
This paper introduces two new distributed algorithms for the problem of
synchronizing orthogonal matrices over graphs. 
Synchronization means that compositions of transformations (multiplications of matrices) over loops in the graph equal the identity (matrix)~ \cite{singer2011angular,bandeira2013cheeger,wang2013exact,bernard_cvpr}.
Thus, ``synchronization'' does not refer to the related concepts of consensus~\cite{sab04} or rendezvous, e.g., attitude synchronization~\cite{thunberg2014distributed}.
 We formulate the problem as a nonlinear least-squares optimization with matrix variables~\cite{helmke2012optimization,absil2009optimization}. 
For symmetric communication topologies we provide an algorithm
with strong convergence guarantees -- the solution converges to 
the optimal solution of a spectral relaxation, which in turn is  known to
produce near-optimal solutions. For graphs that are possibly asymmetric we provide an algorithm with weaker convergence guarantees but with good performance in numerical simulations.

The synchronization problem appears as one of the key components in the following applications: the 3D-localization problem, where the transformations are obtained from camera measurements; the generalized Procrustes problem, where scales, rotations, and translations are 
calculated between multiple point clouds~\cite{Gower:2004uu}; the 
image registration problem, where transformations are calculated between multiple 
images~\cite{bernard_sami}. Due to sensor and communication limitations, there is often a need to use distributed protocols for the 3D-localization problem and several approaches have been proposed recently~\cite{aragues2015parallel,tron2014distributed,montijano2014distributed}. There are also many other interesting applications for the synchronization problem,
see Section 1.2 in \cite{boumal2015riemannian}.

If we exclude the requirement that the synchronization method
shall be distributed, there is an extensive body of work.
Govindu et al. have presented several approaches based on Lie-group averaging,
where a first-order approximation in the tangent space is used~\cite{govindu2004lie,govindu2006robustness,govindu2014averaging}.
Singer {et al.}  have presented several optimization approaches~\cite{singer2011angular,bandeira2013cheeger,wang2013exact,Singer:2011ba,Hadani:2011hb,Hadani:2011tw,Chaudhury:2013un}.
Pachauri et al. have addressed the special case where the matrices are permutation 
matrices~\cite{Pachauri:2013wx}. In~\cite{wang2013exact}, three types of relaxations of the problem are presented: semidefinite programming relaxation (see \cite{boumal2015riemannian} for an extensive analysis of this approach); spectral relaxation; least unsquared deviation in combination with semidefinite relaxation.
These three relaxations were evaluated in the probabilistic framework where the error to the ground truth was calculated in numerical experiments. The simulations showed that the first two approaches were on par, whereas the last approach performed slightly better. Furthermore, the last approach was significantly more robust to outliers.
The first distributed algorithm we present has a connection to the second category of the three relaxations above, since 
the matrices in the algorithm converge to the optimal solution of the spectral relaxation. Our methods are extrinsic, in the sense that the matrices are 
calculated in $\mathbb{R}^{d \times d}$ and then projected onto the 
set of orthogonal matrices. The opposite to extrinsic methods are intrinsic methods where no projections from an ambient space occur. In \cite{afsari2013convergence}, intrinsic 
gradient descent methods are studied for the problem of finding the 
Riemannian center of mass. 

The contributions of this work can be summarised as the introduction of two novel algorithms (Algorithm 1 and 2) for \emph{distributed} synchronization of orthogonal matrices over \emph{directed} graphs.
For both algorithms we provide conditions for guaranteed convergence. The main result of the paper is the above-mentioned convergence in Algorithm 1 to the optimal solution of the spectral relaxation problem (Proposition~\ref{prop:r1:4}).
Previous works in the context of distributed algorithms have focused on undirected graphs and 3D rotations~\cite{aragues2015parallel,tron2014distributed,montijano2014distributed}.
However, in this work we consider directed graphs and arbitrary dimensions. It should be noted that some of the existing algorithms can be extended to higher dimensions and are given for the 3D-case mostly for clarity of exposition. 

The distributed approaches in this work
bear a resemblance to linear consensus protocols~\cite{mesbahi2010graph,jad03,olfati2007consensus,olfati2004consensus}. 
The methods also share similarities with
the eigenvector method in~\cite{howard2010estimation} and gossip 
algorithms~\cite{boyd2006randomized}.
The important states in our algorithms are matrices, and those combined converge
to a tall matrix whose range space is a certain linear subspace. In the case of symmetric communication between agents, the proposed method can either be interpreted as an extension of the power method or the steepest descent method. In our methods, instead of using the graph Laplacian 
matrix~\cite{mesbahi2010graph}, matrices similar to the 
graph connection Laplacian matrix~\cite{singer2012vector} are used. These matrices can be
seen as a generalizations of the graph Laplacian matrix, in which the scalars 
are replaced by matrix blocks.

The paper proceeds as follows. In Section \ref{sec:2} we introduce the definitions that are necessary in order to precisely state the problem, which is done in Section \ref{sec:problem}. Subsequently, the distributed method for the case of symmetric graphs (Algorithm 1) is introduced and analyzed in section \ref{sec:symmetric}. In Section \ref{sec:asymmetric}, the distributed method for the case of directed and possibly asymmetric graphs (Algorithm 2) is introduced and analyzed. In Section~\ref{sec:conclusion}, the paper is concluded.

\section{Preliminaries}\label{sec:2}
\subsection{Directed Graphs}
Let $\mathcal{G} = (\mathcal{V}, \mathcal{E})$ be a directed graph,
where $\mathcal{V} = \{1, 2, \ldots, n\}$ is the node set and 
$\mathcal{E} \subset \mathcal{V} \times \mathcal{V}$ is the 
edge set. Throughout the paper, the notation $\mathcal{A} \subset \mathcal{B}$ means that every element in $\mathcal{A}$ is contained in $\mathcal{B}$. The set $\mathcal{N}_i$ is the set of neighboring nodes of node $i$ and	 defined by 
\begin{equation}
\mathcal{N}_i = \{j: (i,j) \in \mathcal{E}\}.
\end{equation}
The adjacency matrix $A = [A_{ij}]$ for the graph 
$\mathcal{G}$ is defined by
\begin{equation}
{A}_{ij}= \begin{cases}
1 &  \text{ if } (i,j) \in \mathcal{E},\\
0 & \text{ else.}
\end{cases}
\end{equation}
The graph Laplacian matrix is defined by 
\begin{equation}
L = \text{diag}(A 1_n) - A,
\end{equation}
where $1_n \in \mathbb{R}^n$ is a vector 
with all entries equal to $1$. In order to emphasize that 
the adjacency matrix $A$, the graph Laplacian matrix $L$ 
and the $\mathcal{N}_i$ sets depend on the graph $\mathcal{G}$, 
we may write $A(\mathcal{G})$, $L(\mathcal{G})$ and $\mathcal{N}_i(\mathcal{G})$ respectively. For simplicity however, we mostly omit this notation
and simply write $A$, $L$, and $\mathcal{N}_i$.

\begin{defn} \emph{(connected graph, {undirected} path)}\\
The directed graph $\mathcal{G}$ is connected 
if there is an {undirected} path from any node
in the graph to any other node. An {undirected} path is defined as a (finite) sequence of unique nodes such that for any pair $(i,j)$ of consecutive nodes in the sequence it holds that
$((i,j) \in \mathcal{E}) \text{ or } ((j,i) \in \mathcal{E}).$
\end{defn}

\begin{defn} \emph{(quasi-strongly connected graph, center, directed path)}\label{qscGraph}\\
The directed graph $\mathcal{G}$ is quasi-strongly connected (QSC)
if it contains a center. A center is
a node in the graph to which there is a directed path from any other node
in the graph. A directed path is defined as a (finite) sequence of unique nodes such that any pair of consecutive nodes in the sequence comprises an edge in 
$\mathcal{E}$.
\end{defn}

\begin{defn} \emph{(strongly connected graph)}\label{scGraph}\\
The directed graph $\mathcal{G}$ is strongly connected
if for all pairs of nodes $(i,j) \in \mathcal{V} \times \mathcal{V}$, there is a directed path from
$i$ to $j$. 
\end{defn}

\begin{defn}\emph{(symmetric graph)}\label{def:3}\\
The directed graph $\mathcal{G} = (\mathcal{V}, \mathcal{E})$ is {symmetric} if
\begin{equation}
((i,j) \in \mathcal{E}) \Rightarrow ((j,i) \in \mathcal{E}) \text{ for all } (i,j) \in \mathcal{V} \times \mathcal{V}.
\end{equation} 
\end{defn}
Given a graph $\mathcal{G} = (\mathcal{V}, \mathcal{E})$, the graph $\bar{\mathcal{G}} = (\mathcal{V}, \bar{\mathcal{E}})$ is the graph constructed by reversing the direction
of the edges in $\mathcal{E}$, i.e., $(i,j) \in \bar{\mathcal{E}}$ if and only if
$(j,i) \in {\mathcal{E}}$. 
It is easy to see that 
\begin{equation}
A(\bar{\mathcal{G}}) = (A({\mathcal{G}}))^T \text{ and } L(\bar{\mathcal{G}}) =  \text{diag}((A({\mathcal{G}}))^T 1_n) - A({\mathcal{G}})^T.
\end{equation}

\subsection{Synchronization or transitive consistency of matrices}
The set of invertible matrices in $\mathbb{R}^{d \times d}$ is 
$GL(d, \mathbb{R})$ and the
group of orthogonal matrices in $\mathbb{R}^{d \times d}$ is
\begin{equation}
O(d) = \{R \in \mathbb{R}^{d \times d} : R^TR = I_d\}.
\end{equation}
The set $SO(d)$ comprises those matrices in $O(d)$ whose determinants are equal to $1$.
\begin{defn} \emph{(transitive consistency)} \label{def:transcons}
\begin{enumerate}
\item The matrices in the collection $\{{R}_{ij}\}_{(i,j) \in \mathcal{V} \times \mathcal{V}}$ of matrices in $GL(d,\mathbb{R})$ are transitively 
consistent for the complete graph if 
\begin{equation}
{R}_{ik} = {R}_{ij}{R}_{jk}
\end{equation}
for all $i,j$ and $k$. \\
\item Given a graph $\mathcal{G} = (\mathcal{V}, \mathcal{E})$, the matrices in the collection 
$\{{R}_{ij}\}_{(i,j) \in \mathcal{E}}$ of matrices in $GL(d,\mathbb{R})$ are  \emph{transitively 
consistent for $\mathcal{G}$} if there is a collection
$\{{R}_{ij}\}_{(i,j) \in \mathcal{V} \times \mathcal{V}} \supset \{{R}_{ij}\}_{(i,j) \in \mathcal{E}}$
such that $\{{R}_{ij}\}_{(i,j) \in \mathcal{V} \times \mathcal{V}}$ is transitively consistent for the complete graph.
\end{enumerate}
\end {defn}

If it is apparent by the context, sometimes we will be less strict and omit to mention which graph a collection of transformations is transitively consistent for. Another word for transitive consistency is synchronization. We will use the two interchangeably.
 A sufficient
condition for synchronization of the $R_{ij}$-matrices for
any graph is that
there is a collection $\{R_i\}_{i \in \mathcal{V}}$ of matrices in $GL(d,\mathbb{R})$ such that 
\begin{equation}\label{eq:alf:2}
R_{ij} = R^{-1}_{i}R_{j}
\end{equation}
for all $(i,j) \in \mathcal{E}$. Lemma~\ref{lem:1} below and the proof thereof provides additional important information. 
The result is similar to that in \cite{tron2014distributed}.
For the statement of the lemma, the following 
definition is needed. 

\begin{defn}
Two collections $\{R_i\}_{i \in \mathcal{V}}$ and 
$\{\bar{R}_{i}\}_{i \in \mathcal{V}}$ of
matrices in $GL(d,\mathbb{R})$ are equal up to transformation 
from the left, if there is $Q \in GL(d,\mathbb{R})$ such that 
\begin{equation}
QR_i = \bar{R}_i \text{ for all  } i.
\end{equation}
\end{defn}

\begin{lem}\label{lem:1}
For any graph $\mathcal{G} = (\mathcal{V}, \mathcal{E})$ and collection $\{{R}_{ij}\}_{(i,j) \in \mathcal{E}}$ of matrices
in $GL(d,\mathbb{R})$ that are transitively consistent for $\mathcal{G}$,
\begin{enumerate}
\item there is a collection $\{R_{i}\}_{i \in \mathcal{V}}$ of matrices in $GL(d,\mathbb{R})$ such that  
\begin{equation}\label{eq:1}
R_{ij} = R^{-1}_{i}R_{j} \text{ for all } (i,j) \in \mathcal{E},
\end{equation}
\item all collections $\{ R_i\}_{i \in \mathcal{V}}$ satisfying \eqref{eq:1} are equal up
to transformation from the left
if and only if $\mathcal{G}$ is connected,
\item \vspace{1mm} there is a unique collection $\{{R}_{ij}\}_{(i,j) \in \mathcal{V} \times \mathcal{V}} \supset \{{R}_{ij}\}_{(i,j) \in \mathcal{E}}$ of
transitively consistent matrices for the 
complete graph, if and only if all collections $\{ R_i\}_{i \in \mathcal{V}}$ satisfying \eqref{eq:1} are equal up
to transformation from the left.

\end{enumerate}
\end{lem}

\noindent\emph{Proof: } See \cite{thunberg2015transitive}.

Another equivalent definition of transitive consistency or synchronization is given in 
\cite{tron2014distributed,carlone2014angular}. A set of transformations is transitively consistent if the composite
transformations equal the identity along loops or cycles in the graph. In Proposition 7 in \cite{tron2014distributed}
the equivalence between this condition for the loops and \eqref{eq:1} is shown. The definition using the auxiliary $R_i$-matrices, \eqref{eq:1}, is the one we will use mostly in our 
analysis. 
\hfill $\blacksquare$\newline

\section{Problem formulation}\label{sec:problem}
The optimization problem of interest is given by
\begin{equation*}
(P_1)  
\begin{cases}
\begin{aligned}
& \underset{R}{\text{minimize}}
& & f_1(R) =  \sum_{(i,j) \in \mathcal{E}}\frac{a_{ij}}{2}\|R_{ij} - R^{T}_iR_j \|_F^2, \quad \quad \quad \\
& \text{subject to}
& & R  = [R_1, R_2, \ldots, R_n],~\\
& & & R_i  \in O(d) \text{ for all } i,
\end{aligned}
\end{cases}
\end{equation*}
where the $a_{ij}$'s are positive scalar weights, the set $\mathcal{E}$ is the edge-set of a connected directed graph $\mathcal{G} = (\mathcal{V}, \mathcal{E})$, and the matrices in the collection $\{R_{ij}\}_{(i,j) \in \mathcal{E}}$ belong to $O(d)$. The objective function $f_1$ comprises the weighted element-wise sum of squared differences between the $R_{ij}$-matrices and the $(R_i^{T}R_j)$-matrices. The problem $(P_1)$ is similar to the problem in \cite{wang2013exact}. The differences are that we allow for directed graphs (instead of undirected graphs) and we do not require the matrices to be contained in $SO(d)$. 

The overall problem addressed in this paper is how to design  distributed methods that achieve good solutions to $(P_1)$. 

When the $R_i$-matrices are orthogonal, $f_1$ can be written as 
\begin{align}
\label{eq:olle:501}
f_1(U_1) & = \frac{1}{2}\text{tr}(U_1^TL_{\text{undir}}U_1), \text{ where }\\
U_1 & = \begin{bmatrix}
R^{-T}_1 & R^{-T}_2 & \ldots & R^{-T}_n 
\end{bmatrix}^T.
\end{align}
The matrix $U_1$ will be used frequently in the following.
The presented definition of $U_1$ might seem overly complicated, since $R_i^{-T} = R_i$
when $R_i$ is orthogonal. However, we will also use $U_1$ when the $R_i$-matrices are not orthogonal. In that case it is important to note however that $f_1$ is given by the definition in $(P_1)$ and not by \eqref{eq:olle:501}.

The matrix $L_{\text{undir}}$ is
defined as 
\begin{align}
\label{eq:alf:1}
L_{\text{undir}}  = \text{diag}(A{1}_n) \otimes I_d
 +\text{diag}(\bar{W}^T\bar{W}) - (W + W^T),
\end{align} 
where \small
$[W]_{ij} = \begin{cases}
a_{ij}R_{ij} & \text{ if } (i,j) \in \mathcal{E}, \\
0 & \text{ else, }
\end{cases}, \: [A]_{ij} = \begin{cases}
a_{ij} & \text{ if } (i,j) \in \mathcal{E}, \\
0 & \text{ else, }
\end{cases}$  
$[\bar W]_{ij} = \begin{cases}
\sqrt{a_{ij}}R_{ij} & \text{ if } (i,j) \in \mathcal{E}, \\
0 & \text{ else, }
\end{cases}$
\normalsize
and the operator $\text{diag}(\cdot)$ in the second term is understood in the block-matrix sense, i.e. $\text{diag}(W^TW)  = (I_n \otimes 1_d 1_d^T ) \odot W^TW$, where $\odot$ denotes element-wise multiplication.
The matrix $A$ is now, compared to Section~\ref{sec:2}, a weighted adjacency matrix of $\mathcal{G}$. In the following, $A$ will always be defined in this way. 
The matrices $L_{\text{undir}}$, $W$,  and $U_1$ are de-facto functions of the graph $\mathcal{G}$, $A$, and either the $R_{ij}$-matrices or the $R_i$-matrices. However, unless it is absolutely necessary, we will not show this dependence explicitly.

\section{Symmetric Graphs}\label{sec:symmetric}
In this section we introduce Algorithm 1. It is the proposed distributed algorithm for synchronization over symmetric graphs. A detailed analysis of Algorithm 1 will be conducted in Section~\ref{sec:r1:alg1}. 

\subsection{The algorithm}
There are four matrices 
that can be seen as the output of the algorithm at each iteration.
Those are: $\tilde{R}_i(k)$, ${R}_i(k)$, $\tilde{Q}_i(k)$, and ${Q}_i(k)$. The procedure to calculate the $\tilde{R}_i(k)$-matrices is similar to a gradient descent procedure and can also be seen as the power method. The ${R}_i(k)$-matrices are the projections of the $\tilde{R}_i^T(k)$-matrices onto $O(d)$. 

For all $i$, the matrix $\tilde{Q}_i(k)$ and the corresponding ${Q}_i(k)$ is calculated from 
auxiliary variables. The most important such auxiliary variables are the $d_{is}(k)$'s, which are 
calculated in a distributed manner. The protocol for calculating the $d_{is}(k)$'s is similar to a well-known average consensus protocol, but differs by an extra term $(\tilde{d}_{is}(k-1) - \tilde{d}_{is}(k-2))$. This extra term makes the states converge not to the averages of the initial conditions, but to the averages over $j$ for the converging sequences $\{\tilde{d}_{js}(k)\}_{k = 1}^{\infty}$. 
The idea behind the 
$\tilde{Q}_i(k)$-matrices is to modify the $\tilde{R}_i(k)$-matrices in such a way that the modified matrices (the $\tilde{Q}_i(k)$'s) converge to the matrices in the optimal solution to a spectral relaxation of problem ($P_1$) (this relaxation is $(P_2)$ defined in Section~\ref{sec:r1:conv}). 
\rule{\columnwidth}{1pt}
\textbf{Algorithm 1} Distributed method for symmetric graphs 
\hrule
\textbf{Inputs:} a symmetric directed graph $\mathcal{G} = (\mathcal{E}, \mathcal{V})$, a weight matrix $A = [a_{ij}]$, and a collection $\{R_{ij}\}_{(i,j) \in \mathcal{E}}$ of matrices in $O(d)$.\\ \\
\textbf{Outputs:}  $\tilde{R}_i(k)$, ${R}_i(k)$, $\tilde{Q}_i(k)$, ${Q}_i(k)$ for $i \in \mathcal{V}$ and $k \geq 1$.\\ \\
\textbf{Initialization:} let $\tilde{R}_i(0) = I_d$, ${d}_{is}(0) = 1$, $\tilde{d}_{is}(0) = 1$, and $\tilde{d}_{is}(-1) = 1$ for all $i,s$. Let  $V_{ij} = (a_{ij} + a_{ji})I$ 
and $Q_{ij} = a_{ij}R_{ij} + a_{ji}R_{ji}^T$ for all $(i,j) \in \mathcal{E}$. Let $\epsilon_1, \epsilon_2 > 0$.\\

\textbf{Iteration $k \geq 1$:} \\
for all $i$, let
\begin{align*}
\tilde{R}_i(k) & = \tilde{R}_i(k-1) \\
& \quad + \epsilon_1\sum_{j \in \mathcal{N}_i}(Q_{ij}\tilde{R}_j(k-1) - V_{ij}\tilde{R}_i(k-1)), \\
R_i^T(k) & = \text{Pr}_{\mathit{O}(d)}(\tilde{R}_i(k)), \\ \\
\tilde{d}_{is}(k) & = \{\text{ calculated in Subroutine 1 }\}, \\ 
d_{is}(k) & = d_{is}(k-1) + (\tilde{d}_{is}(k-1) - \tilde{d}_{is}(k-2)) \\
& ~ \hspace{3mm}+ \epsilon_2\sum_{l \in \mathcal{N}_i}(d_{ls}(k-1) - d_{is}(k-1)) \\
&~\hspace{3mm} \text{ for } s = 1,2, \ldots, d, \\
D_i(k) & = \text{diag}(d_{i1}(k), d_{i2}(k), \ldots, d_{id}(k)), \\ \\
\tilde{Q}_i(k) & = \{\text{ calculated in Subroutine 1 }\}, \\
Q_i^T(k) & = \text{Pr}_{\mathit{O}(d)}(\tilde{Q}_i(k)(D_i(k))^{-\frac{1}{2}}),
\end{align*}
where $\text{Pr}_{O(d)}$ is the projection operator (the least squares projection onto $O(d)$
computed by means of Singular Value Decomposition (SVD)). \\
\rule{\columnwidth}{1pt}

The $d_{is}$-variables provide a way of creating the $\tilde{Q}_i(k)$-matrices by re-scaling the columns of the $\tilde{R}_i(k)$-matrices. This re-scaling is necessary to obtain the desired convergence. The ${Q}_i(k)$-matrices are projections onto $O(d)$ of scaled versions of the $\tilde{Q}_i^T(k)$-matrices. Under general conditions, those converge to the projections of the matrices in the optimal solution to the spectral relaxation of problem ($P_1$).

\rule{\columnwidth}{1pt}
\textbf{Subroutine 1} Calculation of $\tilde{Q}_i(k)$ and $\tilde{d}_{is}(k)$ 
\hrule

\textbf{Inputs:} $\tilde{R}_i(k)$. \\ \\
\textbf{Outputs:} $[\tilde{d}_{i1}(k), \tilde{d}_{i2}(k), \ldots, \tilde{d}_{in}(k)]^T$, $\tilde{Q}_i(k)$. \\
\begin{enumerate}
\item \textbf{If} $\tilde{R}_i(k)$ is not invertible or $\tilde{R}_i(k)$ is invertible and it does not hold that the eigenvalues of $\tilde{R}_i(k)^{-1}\tilde{R}_i(k-1)$ are distinct, real, and positive. Let $\tilde{Q}_i(k) = \tilde{R}_i(k)$ and $\tilde{d}_{is}(k) = 1$ for all $s$.  \\
\item \textbf{Else} using eigenvalue decomposition, compute 
\begin{equation}
\tilde{R}_i(k-1)^{-1}\tilde{R}_i(k) = P_i^{-1}(k)\bar{D}_i(k)P_i(k).
\end{equation}
Let $\tilde{Q}_i(k) = \tilde{R}_i(k)P_i^{-1}(k)( \bar D_i(k))^{-k}$.

\noindent Let $\tilde{d}_{is}(k) = \|\tilde{q}_{is}(k)\|^2_2 \text{ for all } s,$
where $\tilde{Q}_{i}(k) = [\tilde{q}_{i1}(k), \tilde{q}_{i2}(k), \ldots, \tilde{q}_{id}(k)]$ and each $\tilde{q}_{is} \in \mathbb{R}^{d}$.
\end{enumerate}
\rule{\columnwidth}{1pt}

Let $\tilde R(k)= [\tilde R_1^T(k),\tilde R_2^T(k), \ldots, \tilde R_n^T(k)]^T$ for all $k$.
The update for $\tilde{R}(k)$ is given by 
\begin{align}
\label{eq:r1:2}
\tilde{R}(k+1) & = (I - \epsilon_1 L_{\text{undir}})\tilde{R}(k).
\end{align}

\subsection{Analysis}\label{sec:r1:alg1}
In this section we show how the matrices $\tilde{R}_i(k)$, ${R}_i(k)$, $\tilde{Q}_i(k)$, and ${Q}_i(k)$ relate to problem $(P_1)$. We will provide conditions for well-posedness and convergence.

\subsubsection{Some properties of $L_{\text{undir}}$}\label{sec:r1:5}
For the analysis of Algorithm 1 we first provide an alternative definition of transitive consistency, formulated in terms of the $L_{\text{undir}}$-matrix. To be more precise, in Proposition~\ref{prop:2} we state that for the general case of invertible matrices, transitive consistency is equivalent to the $L_{\text{undir}}$-matrix having a $d$-dimensional nullspace. In other words, there are no collections of matrices that are not transitively consistent for which the $L_{\text{undir}}$-matrix has a nullspace of dimension $d$. This motivates  the choice of $f_1$ as the objective function in an optimization problem for synchronization of matrices.

\begin{prop}\label{prop:2}
For collections $\{{R}_{ij}\}_{(i,j) \in \mathcal{E}}$ of matrices in $GL(d,\mathbb{R})$ and graph $\mathcal{G} = (\mathcal{V}, \mathcal{E})$ that is connected, it holds that 
\begin{equation}
{\text{dim}(\text{ker}(L_{\text{undir}}})) \leq d.
\end{equation}
with equality if and only if transitive consistency holds.
\end{prop}

Before we provide the proof of Proposition~\ref{prop:2}
we provide the following lemma and the proof thereof.
\begin{lem}\label{lem:5}
For any connected graph $\mathcal{G} = (\mathcal{V}, \mathcal{E})$
and collection $\{{R}_{ij}\}_{(i,j) \in \mathcal{E}}$ of matrices
in $GL(d, \mathbb{R})$, the collection $\{{R}_{ij}\}_{(i,j) \in \mathcal{E}}$ is transitively consistent
for $\mathcal{G}$
if and only if there is a collection $\{R_i\}_{i\in \mathcal{V}}$ of matrices in $GL(d, \mathbb{R})$ such that
\begin{equation} \label{eq:olle:500}
\text{im}(U_1) \subset \text{ker}(L_{\text{undir}}).
\end{equation}
\end{lem}

\noindent\emph{Proof: }
Suppose $\{{R}_{ij}\}_{(i,j) \in \mathcal{E}}$ is transitively consistent, then, according to Lemma~\ref{lem:1}, there is  $\{R_i\}_{i\in \mathcal{V}}$ such that \eqref{eq:1} holds for the $R_i$-matrices.  In this case it holds that 
\begin{align}
\label{eq:3}
U_1^TL_{\text{undir}}U_1 = 0,
\end{align}
which implies that \eqref{eq:olle:500} is fulfilled since $L_{\text{undir}}$ is symmetric.
On the other hand, if $\{R_{ij}\}_{(i,j) \in \mathcal{E}}$ is not transitively 
consistent, there are no $R_{i}$ such that \eqref{eq:1} holds. It can now be shown that \eqref{eq:3} does not hold for any collection $\{R_{i}\}_{i \in \mathcal{V}}$ of matrices in $GL(d, \mathbb{R})$.
\hfill $\blacksquare$\newline

\noindent\emph{Proof of Proposition~\ref{prop:2}: } \\
$\textbf{Part 1:}$
Here we assume that $\{{R}_{ij}\}_{(i,j) \in \mathcal{E}}$ is 
transitively consistent.
Due to Lemma~\ref{lem:5}, we know that 
\begin{equation}\label{eq:4}
\text{dim}(\text{ker}(L_{\text{undir}})) \geq d.\end{equation}
Thus we need to show that the inequality in \eqref{eq:4} cannot be strict.
Since $\{{R}_{ij}\}_{(i,j) \in \mathcal{E}}$ is transitively consistent, there is $\{R_i\}_{i\in \mathcal{V}}$, where the $R_i$ fulfill \eqref{eq:1}.

Suppose the inequality \eqref{eq:4} is strict. 
We know that
$\text{im}(U_1) \subset \text{ker}(L_{\text{undir}}).$
Now there must be a vector $y = [y_1^T, y_2^T, \ldots, y_n^T]^T \in \mathbb{R}^{nd}$, where the $y_i$ are in $\mathbb{R}^d$, such that 
 $y \in \text{ker}(L_{\text{undir}}),$ $y \neq 0$, and $y^TU_1 = 0$. 
 There must be $k$ and $l$ such that the $l$-th element of $y_k$ is nonzero. 
 Now, 
 let  
 \begin{align*}
\bar {X} & = [\bar{x}_1, \bar{x}_2, \ldots, \bar{x}_d] \\
& = [(R^{-1}_kR_1)^{-T}, (R^{-1}_kR_2)^{-T}, \ldots, (R^{-1}_kR_n)^{-T}]^T
 \end{align*}
 and 
  $\bar{Y} = [\bar{x}_1, \bar{x}_2, \ldots, \bar{x}_{l-1}, y, \bar{x}_{l+1}, \bar{x}_d],$
where $\bar{x}_s \in \mathbb{R}^{nd}$ for $s = 1, 2, \ldots n$.
It holds that $y^T\bar X = 0$ and 
$L_{\text{undir}}\bar{Y} = 0$. For all $i$, let $\bar{Y}_i$
denote the $i$-th $d \times d$ block matrix in $\bar{Y}$. 
The rest of this part of the proof consists of firstly showing that all the $\bar Y_i$-matrices are invertible and secondly showing that we can use those matrices to formulate a contradictory statement.

It holds that
$\bar{Y}_{k} \in GL(d, \mathbb{R})$. This is true since it is constructed by taking the identity matrix and replacing the $l$-th column by another vector that has a nonzero $l$-th element.
Now, for any $j \in \mathcal{N}_k$
it holds that 
${\|R_{kj}\bar{Y}_j - \bar{Y}_k\|_F = 0,}$
which implies that $\bar{Y}_j \in GL(d, \mathbb{R})$.
Also, for any $i$ such that $k \in \mathcal{N}_i$, it holds that 
$\|R_{ik}\bar{Y}_k - \bar{Y}_i\|_F = 0,$
which implies that $\bar{Y}_i \in GL(d, \mathbb{R})$. Now, due to the fact that $\mathcal{G}$ is connected, an induction argument
can be used to show that all the $\bar{Y}_i$ are elements in $GL(d, \mathbb{R})$.

The collection $\{\bar{Y}_i\}_{i \in \mathcal{V}}$ satisfies 
$R_{ij} = \bar{Y}_{i}\bar{Y}_{j}^{-1}$ for all $(i,j) \in \mathcal{E}.$
Since $y \not\in \text{im}(U_1)$, the two collections $\{\bar{Y}^{-1}_i\}_{i \in \mathcal{V}}$ 
and $\{R_i\}_{i \in \mathcal{V}}$ are not equal up to transformation from the left. But, since the graph is connected, the two must be equal up to 
transformation from the left (Lemma \ref{lem:1}). This is a contradiction.  
Hence it is a false assumption that the inequality in \eqref{eq:4} is strict.

$\textbf{Part 2:}$ Here we show that if ${\text{dim}(\text{ker}(L_{\text{undir}}}))  \geq d,$ then $\{R_{ij}\}_{(i,j) \in \mathcal{E}}$ is transitively 
consistent. 

Let $\bar{Y} = [y_1, y_2, \ldots, y_{nd}]^T \in \mathbb{R}^{nd \times d}$ be any full rank matrix such that 
$L_{\text{undir}}\bar{Y} = 0.$
It holds that all the $y_i \in \mathbb{R}^{d}$.
Let $\bar{Y}_i$
be the $i$-th $d \times d$ block matrix in $\bar{Y}$. 
Since $\bar{Y}$ is full rank, there is a collection  
$\{y_{i_j}\}_{j =1}^d$ such that $[y_{i_1}, y_{i_2}, \ldots
y_{i_d}] \in GL(d, \mathbb{R})$. 

Now, for $k \in \mathcal{V}$ we know that for any $j \in \mathcal{N}_k$
it holds that 
$\|R_{kj}\bar{Y}_j - \bar{Y}_k\|_F = 0,$
which implies that $\text{im}(\bar{Y}_j^T) = \text{im}(\bar{Y}_k^T)$.
Also, for any $i$ such that $k \in \mathcal{N}_i$, it holds that 
$$\|R_{ik}\bar{Y}_k - \bar{Y}_i\|_F = 0,$$
which implies that $\text{im}(\bar{Y}_i^T) = \text{im}(\bar{Y}_k^T)$. Now, due to the fact that $\mathcal{G}$ is connected, an induction argument
can be used to show that $\text{im}(\bar{Y}_j^T) = 
\text{im}(\bar{Y}_i^T)$ for all $i,j$. But then 
$\text{im}([y_{i_1}, y_{i_2}, \ldots
y_{i_d}]) \subset \text{im}(\bar{Y}_j^T) \text{ for all }j,$
which together with the fact that $[y_{i_1}, y_{i_2}, \ldots
y_{i_d}] \in GL(d, \mathbb{R})$ is full rank, 
implies that $\bar{Y}_i \in GL(d, \mathbb{R})$ for all $i$.
It holds that
$\text{im}(U_1(\{\bar{Y}^{-1}_i\}_{i \in \mathcal{V}})) \subset  \text{ker}(L_{\text{undir}}).$
Now the desired result follows by application of Lemma~\ref{lem:5}.
\hfill $\blacksquare$\newline

From Definition \ref{def:transcons} and Proposition \ref{prop:2} we get the following equivalent characterisations of transitive consistency. 
\begin{cor}{(equivalent characterisations of transitive consistency)}\label{cor:r1:1}\\
For a connected graph $\mathcal{G} = (\mathcal{V}, \mathcal{E})$ and a collection $\{{R}_{ij}\}_{(i,j) \in \mathcal{E}}$ of matrices in $GL(d,\mathbb{R})$ the following three statements are equivalent
\begin{enumerate}
\item $\{{R}_{ij}\}_{(i,j) \in \mathcal{E}}$ is transitively consistent. \\
\item There is a collection $\{\bar R_i\}_{i \in \mathcal{V}}$ of matrices in $GL(d, \mathbb{R})$ such that $\text{ker}(L_{\text{undir}}(\{{R}_{ij}\}_{(i,j) \in \mathcal{E}})) = \text{im}(U_1(\{\bar R_i\}_{i \in \mathcal{V}}))$. \\
\item There is a collection $\{\bar R_i\}_{i \in \mathcal{V}}$ of matrices in $GL(d, \mathbb{R})$ such that $R_{ij} = \bar R^{-1}_{i}\bar R_{j} \text{ for all } (i,j) \in \mathcal{E}.$
\end{enumerate} 
\end{cor}
According to Corollary~\ref{cor:r1:1}, the following holds.
For a collection $\{\bar R_i\}_{i \in \mathcal{V}}$ of matrices in $GL(d, \mathbb{R})$, $f_1$ attains the value $0$ if and only if the collection $\{R_{ij}\}_{(i,j) \in \mathcal{E}}$ is transitively consistent. This means that minimization of the right-hand side of \eqref{eq:olle:501} is an approach to consider even when the $R_{ij}$ are not necessarily orthogonal. This is the approach in the first step of an iterative method recently published~\cite{Johan_mtns_2016_2}.
   
With the assurance given by Proposition~\ref{prop:2} that $f_1$ is a suitable objective function, we now move on to the convergence analysis of Algorithm 1.

\subsubsection{Convergence analysis}\label{sec:r1:conv}
We begin by introducing
a relaxation of  
problem $(P_1)$, given by
\begin{equation*}
(P_2)  
\begin{cases}
\begin{aligned}
& \underset{X \in \mathbb{R}^{nd \times d}}{\text{minimize}}
& & \frac{1}{2}\text{tr}(X^TL_{\text{undir}}X)\\
& \text{subject to}
& & X^TX = nI_d.  \quad \quad \quad \quad  \quad \quad \quad \quad  \quad \quad \quad \quad  \quad \quad \quad 
\end{aligned}
\end{cases}
\end{equation*}
Let $\bar{X} =[\bar R_1^T, \bar R_2^T, \ldots, \bar R_n^T]^T$ be one of the optimal solutions to $(P_2)$.
Under the assumption that $R_{ij} = R_{ji}^T$ for all $(i,j) \in \mathcal{E}$, the spectral relaxation method in Section 2.2 in \cite{wang2013exact} is the same as solving problem $(P_2)$.

Now we provide a list of conditions for convergence, which are recalled in the following propositions. Only a subset of the conditions will be used in each proposition.

\begin{table}[!ht]
  \begin{center}
    \begin{tabular}{| l | l |}
    \hline
    (1) & $\mathcal{G}= (\mathcal{V}, \mathcal{E})$ is connected 
and symmetric. \\
    \hline
    (2) & $R_{ij} \in O(d)$ for all $(i,j) \in \mathcal{E}$.\\
        \hline
    (3) & $\{R_{ij}\}_{(i,j) \in \mathcal{E}}$ is transitively consistent. \\
        \hline
    (4) & $\epsilon_1 < \frac{2}{\|P\|_2}$, where \\
    & $P = \text{diag}((A+A^T)1_n) + A + A^T$.\\
        \hline
    (5) & $\epsilon_1 < \frac{2}{\|\bar L\|_2}$, where \\
    & $\bar L = \text{diag}((A+A^T)1_n) - (A+A^T)$.\\
        \hline
    (6) & $\bar{R}_i \in GL(d, \mathbb{R})$ for all the $\bar{R}_i$-matrices in $\bar X$. \\
        \hline
    (7) & $\sum_{i \in \mathcal{V}} \bar R_i \in GL(d, \mathbb{R}),$
where the $\bar{R}_i$ are the \\
&  matrices in $\bar X$.\\
          \hline
    (8) & It holds that 
$\lambda_{(n-1)d} > \lambda_{(n-1)d+1},$ \\
      &  where $\lambda_1 \geq \lambda_2 \geq \cdots \geq \lambda_n$ are the
eigenvalues \\  
      & of $L_{\text{undir}}$, i.e., 
$L_{\text{undir}} = V\text{diag}([\lambda_1, \lambda_2, \ldots, \lambda_{nd}])V^T$. \\
    \hline
    (9) & It holds that 
$\lambda_{nd-(i+1)} > \lambda_{nd - i},$
for \\ 
      & $i = 0,1,\ldots d-2$  where the $\lambda_i$ are 
\\ 
& defined in (8) above. \\ 
    \hline
   (10) & $\epsilon_2 < \frac{2}{\|L\|_2}$, where $L$ is the graph Laplacian \\ & matrix of the graph $\mathcal{G}$ defined in Section~\ref{sec:2} \\ &(not to mix up with $\bar L$). \\
    \hline
    \end{tabular}
  \end{center}
  \caption{Conditions for convergence.}
\end{table}

The conditions (1-3) are fundamental properties that need 
no further explanation. Conditions (4-5) and (10) are 
conditions for the step size determination. These have the property that they scale with the number of nodes in the network, i.e., $n$.
Condition (6) states that all the $\bar{R}_i$-matrices in the 
optimal solution to $(P_2)$ are invertible and condition (7) states that the sum of those is invertible. Condition (8) states that the $d$ smallest eigenvalues of $L_{\text{undir}}$ are strictly smaller than the $(n-1)d$ largest. Condition (9) states that the $d$ smallest eigenvalues of $L_{\text{undir}}$ are distinct. 
 
The following lemma provides a bound 
for $\epsilon_1$ such that the discrete-time system defined by \eqref{eq:r1:2} is stable. It is a justification of convergence condition (4).

\begin{lem}\label{lem:29}
The largest eigenvalue of $P$, as defined in convergence condition (4), is an upper bound for the eigenvalues 
of $L_{\text{undir}}$ for all graphs $\mathcal{G} = (\mathcal{V}, \mathcal{E})$ and collections $\{R_{ij}\}_{(i,j) \in \mathcal{E}}$ satisfying convergence conditions (1) and (2).
\end{lem}

\noindent\emph{Proof: } 
The largest eigenvalue $\lambda_{\max}$ 
is given by $$\lambda_{\max} = \max_{v \in \mathbb{S}^{nd}}v^TL_{\text{undir}}v,$$
where $\mathbb{S}^{nd}$ is the $nd$-dimensional unit sphere.
Let $v = [v_1^T, v_2^T, \ldots, v_n^T]^T \in \mathbb{S}^{nd}$, where
each $v_i \in \mathbb{R}^{d}$. By using the structure of the $f_1$-function in $(P_1)$,
one can show that \\
$v^TL_{\text{undir}}v  = \sum\limits_{(i,j) \in \mathcal{E}}a_{ij}\|R_{ij}v_j - v_i \|_2^2.$
Now, 
$v^TL_{\text{undir}}v  \leq \sum\limits_{(i,j) \in \mathcal{E}}a_{ij}(\|v_i\|_2 + \|v_j\|_2)^2
\leq \max_{u\in\mathbb{S}^n \cap (\mathbb{R}^{+})^n}u^TPu.$ The set $\mathbb{R}^+$ comprises the non-negative real numbers.
\hfill $\blacksquare$\newline

Lemma~\ref{lem:29} has the following implication. If $\epsilon_1$ in 
Algorithm 1 is chosen to be smaller than $2/\|P\|_2$, then $\tilde R(k)$ (see \eqref{eq:r1:2}) converges
as $k$ goes to infinity. 

Now, unless the $R_{ij}$-matrices are transitively consistent, the nullspace of $L_{\text{undir}}$ has a lower dimension than $d$ and in general it will be zero-dimensional.
Thus, $\tilde R(k)$ converges to zero. In the case when the $R_{ij}$-matrices are transitively consistent, the $\tilde R(k)$ converge to a $d$-dimensional 
subspace in general.

Now we provide a result for the special case when 
transitive consistency holds. Algorithm 1 reduces to the first two lines in each iteration, which is de facto the power method. See \cite{howard2010estimation} for a discussion about the power method in a similar context. We provide Proposition~\ref{prop:r1:2} and its proof below for the sake of completeness.

\begin{prop}\label{prop:r1:2}
Suppose convergence conditions (1-3), and (5-7) hold.
Then, for Algorithm 1, there is a positive integer $K$ such that $\tilde{R}_j^{-1}(k)$ is well defined for $k\geq K$, and (for $k \geq K$) it holds that
\begin{align*}
(R_i^T(k)R_j(k) \rightarrow R_{ij} \text{ as } k \rightarrow \infty) \text{ for all } (i,j) \in \mathcal{E}, \\
(\tilde{R}_i(k)\tilde{R}_j^{-1}(k) \rightarrow R_{ij} \text{ as } k \rightarrow \infty) \text{ for all } (i,j) \in \mathcal{E}.
\end{align*}
\end{prop}

\noindent\emph{Proof: } 
Under conditions (2) and (3), it holds that $L_{\text{undir}}$ is similar to the 
matrix $\bar{L} \otimes I_d$, where $\bar{L}$ is defined in convergence condition (5). This is a consequence of the fact that we can write $L_{\text{undir}}$ as
$$L_{\text{undir}} = \text{diag}(U_1(\{R_i\}_{i \in \mathcal{V}}))(\bar{L} \otimes I_d)\text{diag}(U_1(\{R_i^T\}_{i \in \mathcal{V}})).$$
In the right-hand side above, $\text{diag}(\cdot)$ shall be interpreted in the 
block diagonal sense, where the $R_i^T$ and $R_i$ are put as blocks on the diagonal.
 
Since conditions (1), (2), and (3) hold, we can use Corollary~\ref{cor:r1:1}. There is a collection $\{\bar R_i\}_{i \in \mathcal{V}}$ such that the conditions \emph{(1)} and \emph{(2)} are fulfilled (in the corollary). Thus,
$\text{ker}(L_{\text{undir}}) = \text{im}(U_1(\{\bar R_i\}_{i \in \mathcal{V}})).$
Due to this fact and the fact that 
$\tilde{R}_i(0) = I_d$ for all $i$, it holds that 
the projection of $\tilde{R}^T(0)$ onto $\text{ker}(L_{\text{undir}})$
is 
$U_1U_1^T \tilde{R}(0)=U_1V,$
where
$V = \frac{1}{n}\sum_{i \in \mathcal{V}}\bar R_i$. Furthermore,
$R(k)$ converges to $U_1(\{\bar R_i\}_{i \in \mathcal{V}})V$ as 
$k$ goes to infinity (condition (5)). Now, if $V$ is invertible (condition (7)), by the definition of the limit there is a $K$ such that $\tilde{R}_j^{-1}(k)$ is well defined for all $k\geq K$. It holds that 
\begin{align}
 \tilde{R}_i(k) & \rightarrow \bar{R}_iV \text{ as } k \rightarrow \infty, \text{ for all } i,\\
 \label{eq:olle:710}
 {R}_i^T(k) & \rightarrow \bar{R}_i\text{Pr}_{O(d)}(V) \text{ as } k \rightarrow \infty \text{ for all } i.
\end{align}
Convergence conditions (6) and (7) holds. Thus, the $(\bar{R}_iV)$-matrices are invertible. 
Under these conditions, the projections of $\tilde{R}_i(k)$ converge in \eqref{eq:olle:710} (see the last paragraph of Proposition~\ref{prop:r1:4} below for details about the convergence of the projections).  
\hfill $\blacksquare$

In Proposition~\ref{prop:r1:2} it is important to note that for all $(i,j) \in \mathcal{E}$, 
$R_i^T(k)R_j(k)$ converges to $R_{ij} \in {O}(d)$ inside $O(d)$, 
whereas $\tilde{R}_i(k)\tilde{R}_j^T(k)$ converges to $R_{ij} \in {O}(d)$ without guarantees of being in $O(d)$ for each $k$. 

Now we take a step further in our analysis of Algorithm 1. We show 
that when the $R_{ij}$-matrices are not necessarily transitively 
consistent, we still have a nice convergence property for the $\tilde{R}_i(k)$-matrices. 

\begin{prop}\label{prop:r1:3}
Suppose that the convergence conditions (1-2), (4), and (6-8) are satisfied.
Then, for Algorithm~1 there is a positive integer $K$ such that $\tilde{R}_j^{-1}(k)$ is well defined for all $k\geq K$, and (for $k \geq K$) it holds that
\begin{equation}
(\tilde{R}_i(k)\tilde{R}_j^{-1}(k) \rightarrow \bar{R}_i\bar{R}_j^{-1}\text{ as } k \rightarrow \infty) \text{ for all } (i,j) \in \mathcal{E}.
\end{equation}
\end{prop}

\noindent\emph{Proof: } Since the convergence conditions (1), (2), and (4) are fulfilled, we know that the discrete time system 
defined in equation \eqref{eq:r1:2} is stable.

The columns of $\bar X$ are, up to scale and orthogonal transformation, the eigenvectors corresponding to the $d$ smallest eigenvalues of $L_{\text{dir}}$.
Now, let us rewrite $\tilde{R}(k)$ as
$\tilde{R}(k) = VV_1(k)V_1 + (1/\sqrt{n})\bar{X}V_2(k)V_2,$
where the columns of the matrix $V$ are the eigenvectors corresponding to the $(n-1)d$ largest eigenvalues of the matrix $L_{\text{undir}}$. Due to convergence condition (8), there are $(n-1)d$ eigenvalues that are strictly larger than the other $d$ eigenvalues.
The matrices $V_1(k)$ and $V_2(k)$ are defined as 
\begin{align}
\label{eq:alf:1}
V_1(k)  =~ & \text{diag}([(1 - \epsilon_1\lambda_1)^k, (1 - \epsilon_1\lambda_2)^k, \\
\nonumber
& \ldots,  (1 - \epsilon_1\lambda_{(n-1)d})^k]), \\
\nonumber
V_2(k)  =~ & \text{diag}([(1 - \epsilon_1\lambda_{(n-1)d+1})^k, (1 - \epsilon_1\lambda_{(n-1)d+2})^k, \\
\label{eq:alf:2}
& \ldots, 
 (1 - \epsilon_1\lambda_{nd})^k]),
\end{align}
and the matrices $V_1$ and $V_2$ are defined as
$V_1 = V^T\tilde{R}(0)$ and $V_2 = \frac{1}{\sqrt{n}}\bar X^T\tilde{R}(0) = \frac{1}{n}\sum_{i \in \mathcal{V}} \bar R_i$. 

Now, by using the definitions of $V_1(k)$ and $V_2(k)$ together with the 
fact that convergence condition (8) holds, we can formulate the following convergence result:
\begin{equation}\label{eq:r1:4}
\rho(k) = \frac{\max\limits_{j}([V_1(k)]_{jj})}{\min\limits_{j}([V_2(k)]_{jj})} \rightarrow 0 \text{ as } k \rightarrow \infty.
\end{equation}
We also see from the definition of $V_2(k)$, that $V_2(k)$ is invertible for all $k$.
For any pair $(i,j) \in \mathcal{E}$, if $\tilde R_j(k)$ is invertible, it holds that
\begin{align}
\label{eq:olle:505}
\tilde{R}_i(k)\tilde{R}_j^{-1}(k) = (\tilde{R}_i(k)V_2^{-1}(k))(\tilde{R}_j(k)V_2^{-1}(k))^{-1}
\end{align} 
and $\tilde R_j(k)$ is invertible if and only if $\tilde{R}_j(k)V_1^{-1}(k)$ is.

For any $i$ it holds that
\begin{align}
\label{eq:olle:506}
&\tilde R_i(k)V_2^{-1}(k) = (1/\sqrt{n})\bar{R}_iV_2(k)V_2V_2(k)^{-1} \\
\nonumber
&  + ([0, \ldots,0, I_d, 0, \ldots, 0] VV_1(k)V^T\tilde{R}(0))V_2(k)^{-1}.
\end{align}
The matrix $(V_2(k)V_2V_2(k)^{-1})$ is similar to $V_2$ for all $k$,
i.e., the eigenvalues are the same for the two. This means that the norm of $((1/\sqrt{n})\bar{R}_iV_2(k)V_2V_2(k)^{-1})$ is bounded from below by some positive constant. 
Now we consider the second term in the sum on the right-hand side of \eqref{eq:olle:506}. Each element in that term is equal to a linear combination of rational expressions
$([V_1(k)]_{jj}/[V_2(k)]_{ll})$ for $j,l \in \{1,2,\ldots,d\}$.
According to \eqref{eq:r1:4}, all those go to zero as $k$ goes to infinity.

Now, if $k$ is large enough, $\rho(k)$ is small enough for all the $\tilde R_i(k)(V_2(k))^{-1}$, and equivalently, $\tilde R_i(k)$ to be invertible. For such $k$, it holds that 
\begin{align*}
~&\tilde{R}_i(k)\tilde{R}_j^{-1}(k) \\
=~ & \bar{R}_i((1/\sqrt{n})V_2(k)V_2V_2(k)^{-1})\cdot \\
  ~& ((1/\sqrt{n})V_2(k)V_2V_2(k)^{-1})^{-1}\bar R_j^{-1} + \mathcal{O}(\rho(k))\\
= ~& \bar{R}_i \bar{R}_j^{-1} + \mathcal{O}(\rho(k))
\end{align*}
for all $(i,j) \in \mathcal{E}$. Thus,
$\tilde{R}_i(k)\tilde{R}_j^{-1}(k) \rightarrow \bar{R}_i \bar{R}_j^{-1}$ as $k \rightarrow \infty$.
\hfill $\blacksquare$\newline

Before we proceed, we introduce $\mathcal{S}_{(P_2)} = \{\bar XQ: Q \in O(d)\}.$ When the convergence conditions (1), (2), and (8) hold, $\mathcal{S}_{(P_2)}$ comprises the entire set of solutions to problem $(P_2)$.

So far we have only addressed the convergence properties for 
$\tilde{R}_i(k)$ and $R_i(k)$. Now we analyze the convergence of $\tilde{Q}_i(k)$ and $Q_i(k)$.
The following result is a main result of this work. In the 
case when the $R_{ij}$-matrices are not necessarily transitively consistent it 
guarantees convergence to the projected matrices obtained from the optimal
solution to problem $(P_2)$. It can be seen as an equivalent result to that of Proposition~\ref{prop:r1:2} for the $R_i(k)$-matrices, but 
this time formulated for the $Q_i(k)$-matrices for the case when transitive consistency does not necessarily hold for the $R_{ij}$-matrices.

\begin{prop}\label{prop:r1:4}
Suppose that the convergence conditions (1-2), (4), and (6-10) are satisfied.
Then, for Algorithm~1 there is a positive integer $K$ such that $\tilde{R}_i^{-1}(k)$ is well defined for all $i$ and $k\geq K$, and (for $k \geq K$) it holds that
\begin{equation}
U_1(\{\tilde{Q}_i(k)\}_{i \in \mathcal{V}}) \rightarrow \mathcal{S}_{(P_2)}\text{ as } k \rightarrow \infty
\end{equation}
and 
\begin{align}
\nonumber
& ({Q}_i^T(k){Q}_j(k) \rightarrow (\textnormal{Pr}_{O(d)}(\bar{R}_i))(\textnormal{Pr}_{O(d)}(\bar{R}_j))^T \text{ as } k \rightarrow \infty) \\
& \text{ for all } (i,j) \in \mathcal{E}.
\end{align}
\end{prop}

\noindent\emph{Proof: } 
The proof will use the notation from Proposition~\ref{prop:r1:3}. According to Proposition~\ref{prop:r1:3} there is a $K$ such that $\tilde{R}_j^{-1}(k)$ is well defined for all $k\geq K-1$. 
In the following we assume that $k \geq K-1$.

$$\tilde{R}_i(k-1)^{-1}\tilde{R}_i(k) = P_i^{-1}(k)\bar D_i(k)P_i(k),$$
where the right-hand side is the Jordan decomposition. 
It holds that 
$$\tilde{R}_i(k-1)^{-1}\tilde{R}_i(k) = \tilde{R}_i(k-1)^{-1}V_2(k)V_2^{-1}(k)\tilde{R}_i(k)$$
and
\begin{align}
\nonumber
&V_2^{-1}(k)\tilde R_i(k) \\
\label{eq:olle:510}
= ~& (1/\sqrt{n})V_2(k)^{-1}\bar{R}_iV_2(k)V_2 \\
\nonumber
&  + V_2(k)^{-1}([0, \ldots,0, I_d, 0, \ldots, 0] VV_1(k)V^T\tilde{R}(0)), \\
\nonumber 
& \\
\nonumber
&V_2^{-1}(k)\tilde R_i(k-1) \\
\label{eq:olle:511}
=~& (1/\sqrt{n})V_2^{-1}(1)V_2(k-1)^{-1}\bar{R}_iV_2(k-1)V_2 \\
\nonumber
&  + V_2^{-1}(1)V_2(k-1)^{-1}\cdot \\
\nonumber
& ([0, \ldots,0, I_d, 0, \ldots, 0] VV_1(k-1)V^T\tilde{R}(0)).
\end{align}
Now, similar to \eqref{eq:olle:506} we see that each right-hand side of \eqref{eq:olle:510} and \eqref{eq:olle:511} a sum of two terms. The first terms have bounded  eigenvalues and the second terms are $\mathcal{O}(\rho(k))$ and $\mathcal{O}(\rho(k-1))$, respectively.
Thus, for all $i$ it holds that 
\begin{align}
\label{eq:0lle:520}
~&\tilde{R}_i^{-1}(k-1)\tilde{R}_i(k) = V_2^{-1}V_2(1)V_2 + \mathcal{O}(\rho(k)),
\end{align}
which means, due to \eqref{eq:r1:4}, that for all $i$ it holds that
\begin{equation}\label{eq:olle:521}
\tilde{R}_i^{-1}(k-1)\tilde{R}_i(k) \rightarrow V_2^{-1}V_2(1)V_2 \text{ as } k \rightarrow \infty.
\end{equation}

Now the following holds:
\begin{align}
\label{eq:olle:600}
\bar {D}_i(k) & \rightarrow V_2(1) \text{ as } k \rightarrow \infty, \\
\label{eq:olle:601}
P_i(k) & \rightarrow \bar{Z}V_2 \text{ as } k \rightarrow \infty,
\end{align}
where $\bar{Z}$ is a diagonal matrix with real nonzero elements on the 
diagonal.
Equation \eqref{eq:olle:600} is a consequence of the continuity of the eigenvalues. For equation \eqref{eq:olle:601} the key point is that
convergence property (9) holds. Thus, the diagonal elements of $V_2(1)$ are distinct and real. Then it holds that for $k$ large enough, the elements of ${D}_i(k)$ are distinct and real. This guarantees the convergence of $P_i(k)$, see 
for example Lemma 12 in~\cite{Anandkumar2012}.  

The matrix $\tilde{Q}_i(k)$ is defined as 
\begin{equation}\label{eq:olle:602}
\tilde{Q}_i(k) = \tilde{R}_i(k)P_i^{-1}(k)( \bar D_i(k))^{-k}.
\end{equation}
It holds that
\begin{align*}
\tilde{Q}_i(k)  =~ & (1/\sqrt{n})\bar{R}_iV_2(k)V_2P_i^{-1}(k)( \bar D_i(k))^{-k} \\
&  + ([0, \ldots,0, I_d, 0, \ldots, 0] VV_1(k)V^T\tilde{R}(0))\\
& \cdot P_i^{-1}(k)(\bar{D}_i(k))^{-k}.
\end{align*}
We will now prove that
\begin{equation}\label{eq:olle:603}
\tilde{Q}_i(k) \rightarrow (1/\sqrt{n})\bar{R}_i\bar{Z}^{-1}\text{ as } k \rightarrow \infty.
\end{equation} 
To prove \eqref{eq:olle:603}, 
we prove that
\begin{equation}\label{eq:polle:1000}
\begin{cases}
V_2(k)V_2P_i^{-1}(k)( \bar D_i(k))^{-k} \rightarrow \bar{Z}^{-1}, \text{ as } k \rightarrow \infty, & \\
VV_1(k)V^T\tilde{R}(0)P_i^{-1}(k)(\bar{D}_i(k))^{-k}  \rightarrow 0 \text{ as } k \rightarrow \infty. &
\end{cases}
\end{equation}
The expressions in
\eqref{eq:polle:1000} 
can be manipulated as follows.
\begin{align*}
& V_2(k)V_2P_i^{-1}(k)( \bar D_i(k))^{-k} \\
=~ & V_2(k)V_2\left (P_i^{-1}(k)( \bar D_i(k))^{-k}P_i(k) \right )P_i^{-1}(k), \\
=~& \left (V_2(k)V_2V_2^{-1}V_2^{-1}(k) \right) \left ( V_2P_i^{-1}(k) \right ) + \mathcal{O}(\rho(k))\\
=~&V_2P_i^{-1}(k) + \mathcal{O}(\rho(k))\\ \\
& VV_1(k)V^T\tilde{R}(0)P_i^{-1}(k)(\bar{D}_i(k))^{-k}\\
=~&  VV_1(k)V^T\tilde{R}(0)\left (P_i^{-1}(k)(\bar{D}_i(k))^{-k}P_i(k) \right )P_i^{-1}(k) \\
=~& \left ( VV_1(k)V^T\tilde{R}(0)V_2^{-1}V_2^{-1}(k)\right ) \left (V_2 P_i^{-1}(k) \right)  + \mathcal{O}(\rho(k))\\
=~& \mathcal{O}(\rho(k)).
\end{align*}
In these derivations we have first used the fact that 
$$(\tilde{R}_i^{-1}(k-1)\tilde{R}_i(k))^{-k} = P_i^{-1}(k)( \bar D_i(k))^{-k}P_i(k).$$
Then we have used \eqref{eq:r1:4}, \eqref{eq:0lle:520}, \eqref{eq:olle:601}, the fact that $V_2(1)^{-k} = V_2^{-1}(k)$, and the fact that the elements of $V_2(k)$ are bounded by $1$ in magnitude.
Now, due to the derivations above we can use \eqref{eq:r1:4} and \eqref{eq:olle:601} to conclude that the convergences in \eqref{eq:polle:1000} hold.

The matrix $D_i(k)$ is defined as 
\begin{equation}\label{eq:olle:604}
D_i(k)  = \text{diag}(d_{i1}(k), d_{i2}(k), \ldots, d_{id}(k)),
\end{equation}
for all $i$.
We will now prove that
\begin{align}
\nonumber
& (d_{is}(k) \rightarrow \frac{1}{n}[\bar{Z}^{-2}]_{ss} \text{ as } k \rightarrow \infty) \\
& \text{ for all } i \in \mathcal{V}, \text{ and } s \in \{1,2, \ldots, d\}.
\end{align}
We recall that $\tilde{d}_{is}(k) = \|\tilde{q}_{is}(k)\|^2_2 \text{ for all } s,$
where $\tilde{Q}_{i}(k) = [\tilde{q}_{i1}(k), \tilde{q}_{i2}(k), \ldots, \tilde{q}_{id}(k)]$ and each $\tilde{q}_{is} \in \mathbb{R}^{d}$.
We see from \eqref{eq:olle:603} that $\tilde{Q}_i(k)$ converges to a constant matrix for all $i$. This means that the $\tilde{d}_{is}$ converge to constants and the averages of those constants are equal to the corresponding elements of $\frac{1}{n}[\bar{Z}^{-2}]_{ss}$. More precisely,
\begin{align}
\label{eq:olle:703}
\frac{1}{n}\sum_{i = 1}^n\tilde{d}_{is}(k) \rightarrow \frac{1}{n^2}[\bar Z^{-1} X^TX \bar Z^{-1}]_{ss} = \frac{1}{n}[\bar Z^{-2}]_{ss}
\end{align}
We will show that for all $i$ and all $s$ it holds that
\begin{align}
\label{eq:olle:605}
\|{d}_{is}(k) - \frac{1}{n}\sum_{j = 1}^n\tilde{d}_{js}(k)\| \rightarrow 0, \text{ as } k \rightarrow \infty,
\end{align}
under the protocol
\begin{align}
\nonumber
d_{is}(k+1) & = d_{is}(k) + (\tilde{d}_{is}(k) - \tilde{d}_{is}(k-1)) \\
\label{eq:olle:611}
& ~ \hspace{3mm}+ \epsilon_2\sum_{l \in \mathcal{N}_i}(d_{ls}(k) - d_{is}(k)) \\
\nonumber
&~\hspace{3mm} \text{ for } s = 1,2, \ldots, d,
\end{align}
and the initial conditions chosen as in Algorithm 1.
The protocol is equal to an average consensus protocol up to the extra term 
$(\tilde{d}_{is}(k) - \tilde{d}_{is}(k-1))$. Without the extra term it is well known that the $d_{is}(k)$'s converge to the averages of the initial conditions under convergence condition (10),
see for example \cite{olfati2007consensus}. With the term and convergence condition (10), the convergence in \eqref{eq:olle:703} holds. Now we show this.

 We can rewrite \eqref{eq:olle:611} as
\begin{align}
\nonumber
& d_{is}(k+1) = d_{is}(k) + \epsilon_2\sum_{l \in \mathcal{N}_i}(d_{ls}(k) - d_{is}(k))  \\ 
\label{eq:olle:610}
& ~ \hspace{1mm}+ \bar{\xi}_s(k) + \bar{\xi}_{is}(k),
\end{align}
where $\bar \xi_s(k) = \frac{1}{n}\sum_{l = 1}^n\xi_{ls}(k)$, $\bar{\xi}_{is}(k) = \frac{1}{n}\sum_{l = 1}^n(\xi_{is}(k) - \xi_{ls}(k))$, and  $\xi_{is}(k) = (\tilde{d}_{is}(k) - \tilde{d}_{is}(k-1))$.
The matrix $L$ is the graph Laplacian matrix of the graph $\mathcal{G}$ (see Section~\ref{sec:2}) and
the matrix $(I_d - \epsilon_2L)$ is a stochastic matrix (convergence condition (10) guarantees that the matrix is stochastic, see \cite{olfati2007consensus}). 
Since we have chosen $d_{is}(0) = 0$, $\tilde{d}_{is}(0) = \tilde{d}_{is}(-1) = 1$ for all $i,s$ it holds that
\begin{align}
\nonumber
& [d_{1s}(k+1), d_{2s}(k+1), \ldots, d_{ns}(k+1)]^T \\
\nonumber
=~& \frac{1}{n}\sum_{i = 1}^n\tilde{d}_{is}(k)[1, 1, \ldots, 1]^T\\
\label{eq:olle:701}
& + \sum_{l = 1}^{k}(I_d - \epsilon_2L)^{k -l}[\bar \xi_{1s}(l), \bar\xi_{2s}(l), \ldots, \bar\xi_{ns}(l)]^T).
\end{align}
The right-hand side of \eqref{eq:olle:701} is written as the sum of two parts. We are done with this section of the proof if it holds that the second part goes to zero as $k$ goes to infinity. Let $P$ be the orthogonal matrix on the left-hand side of in the spectral factorization of $(I_d - \epsilon_2L)$. We can write $P = [P_1, P_2]$, where $P_1$ is a matrix containing the eigenvectors corresponding to eigenvalues in the interior of the unit disc and $P_2$ is the vector corresponding to the eigenvalue that is $1$. We know that such a decomposition exists, see convergence condition (1). It holds that $P_2^T[\bar\xi_{1s}(k),\bar \xi_{2s}(k), \ldots,\bar \xi_{ns}(k)]^T = 0$ for all $k$ and $\|[\bar\xi_{1s}(k), \bar\xi_{2s}(k), \ldots,\bar \xi_{ns}(k)]^T\|$ goes to zero as $k$ goes to infinity.  Thus, after a coordinate transformation with $P$, the second part on the right-hand side of \eqref{eq:olle:701} is equivalent to $n$ convolutions between the impulse responses of stable systems and bounded inputs that converge to zero in norm. Hence the second part goes to zero as $k$ goes to infinity.

Now, by using \eqref{eq:olle:603}, \eqref{eq:olle:703}, and \eqref{eq:olle:605} we deduce that
\begin{equation}\label{eq:olle:606}
\tilde{Q}_i(k)( D_i(k))^{-\frac{1}{2}} \rightarrow \bar{R}_i\bar{P} \text{ as } k \rightarrow \infty,
\end{equation}
where $\bar{P}$ is an orthogonal matrix. By definition it holds that  
\begin{align}
\label{eq:olle:607}
Q_i^T(k)  =~& \text{Pr}_{\mathit{O}(d)}(\tilde{Q}_i(k) D_i^{-\frac{1}{2}}(k)).
\end{align}
$\text{Pr}_{{O}(d)}$ is defined by  $\text{Pr}_{\mathit{O}(d)}(W) =  W_1W_2^T$, where $(W_1,S,W_2)$ is the Singular Value Decomposition (SVD) of $W \in \mathbb{R}^{d \times d}$. We refer to $W_1$ as the ``left matrix'' and $W_2$ as the ``right matrix''. The matrix $\bar R_i\bar P$ is invertible (convergence condition (6)), and for $k$ large enough the sign of the determinant of $(\tilde{Q}_i(k) D_i^{-\frac{1}{2}}(k))$ equals that of $\bar R_i\bar P$. 
This means that for $k$ large enough the projections will be either only onto $SO(d)$ or either only onto $O(d) - SO(d)$. Up to sign the left respective right matrices in the SVD of $(\tilde{Q}_i(k) D_i^{-\frac{1}{2}}(k))$ converge to the left respective the right matrix in the SVD of $\bar R_i\bar P$. When projection is performed, these possible sign differences cancel out in the product between the left and the right matrices.
Thus in \eqref{eq:olle:607}, the expression on the right-hand side goes to $\text{Pr}_{\mathcal{O}(d)}(\bar{R}_i\bar{P}) = \text{Pr}_{\mathcal{O}(d)}(\bar{R}_i)\bar{P}$ $\text{as } k \rightarrow \infty.$ But since $\bar{P}$ is orthogonal, it holds that
\begin{align}
\nonumber
& ({Q}_i^T(k){Q}_j(k) \rightarrow (\textnormal{Pr}_{O(d)}(\bar{R}_i))(\textnormal{Pr}_{O(d)}(\bar{R}_j))^T \text{ as } k \rightarrow \infty) \\
& \text{ for all } (i,j) \in \mathcal{E}.
\end{align}
\hfill $\blacksquare$\newline
The procedure for the $d_{ij}$'s, i.e., to distributed computation of a time-varying signal is equivalent to that in~\cite{zhu2010discrete}.

\subsection{Numerical simulations}\label{sec:num:1}

\begin{figure}[!th] 
\centering
\includegraphics[scale=0.23]{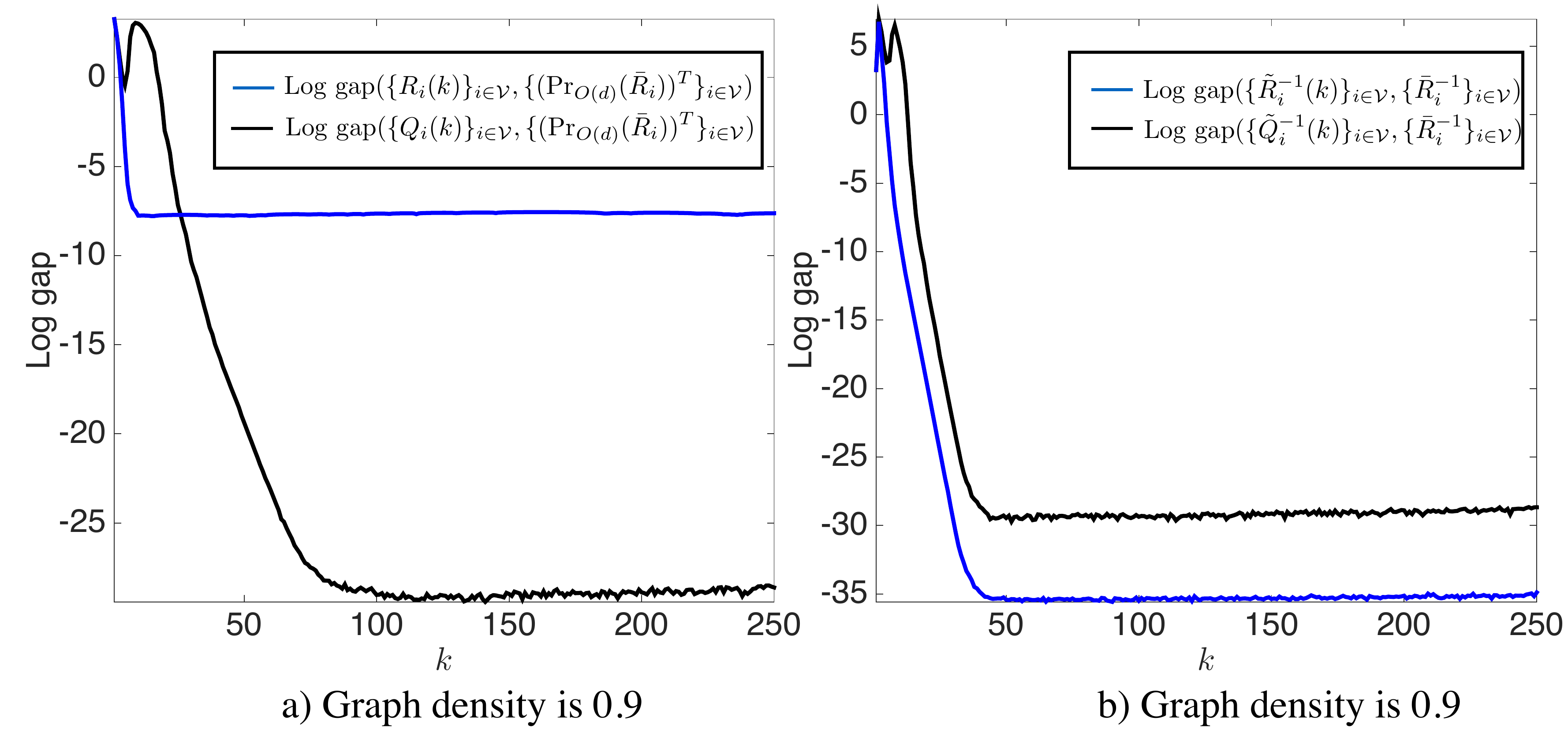}
\centering
\includegraphics[scale=0.23]{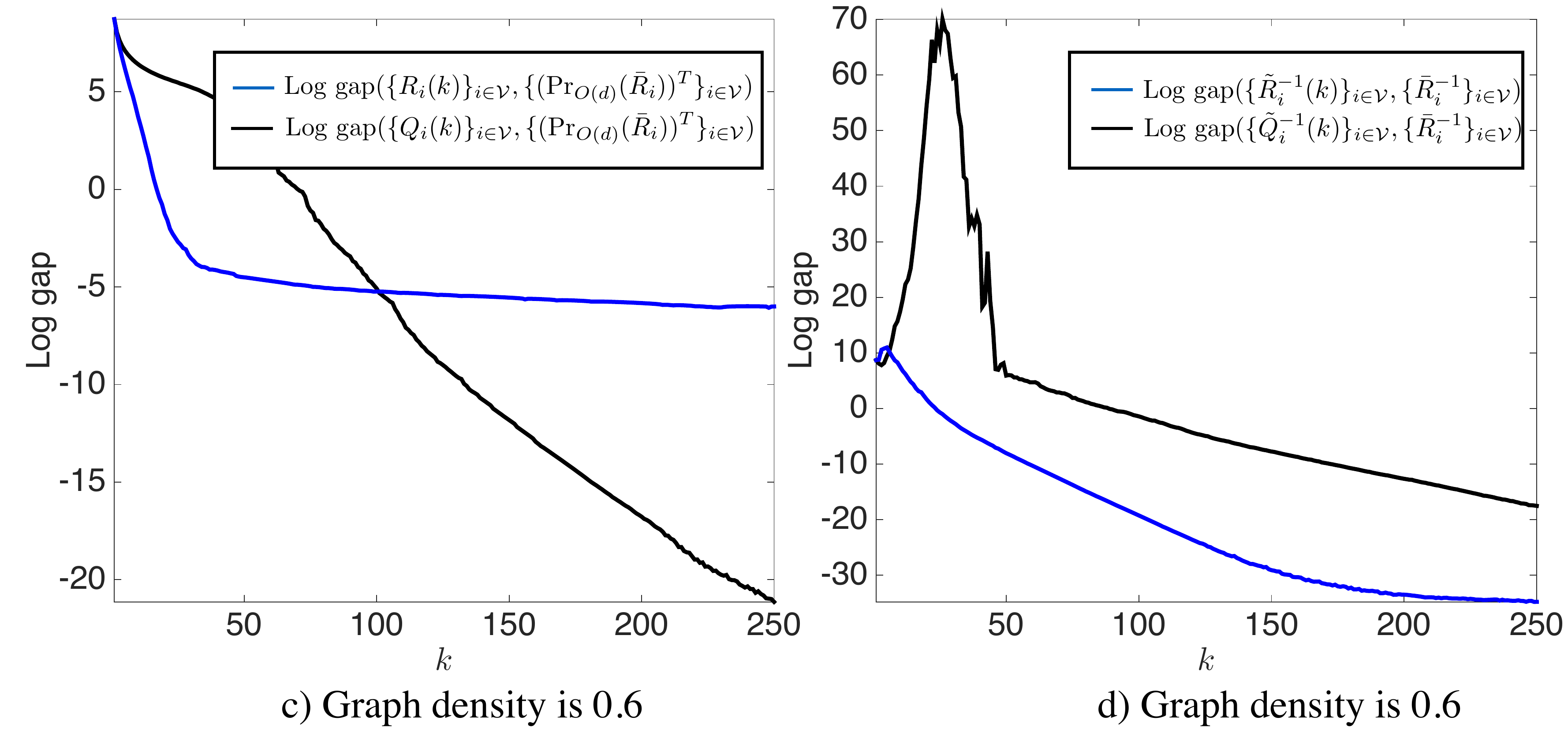}
\caption{Convergence plots for Algorithm 1.}
 \label{fig:1}
\end{figure}

For the collections $\{R_i\}_{i \in \mathcal{V}}$ and $\{\tilde R_i^{-1}\}_{i \in \mathcal{V}}$
of matrices in $GL(d,\mathbb{R})$ we define 
\begin{equation}
\text{gap}(\{R_i\}_{i \in \mathcal{V}}, \{\tilde R_i^{-1}\}_{i \in \mathcal{V}}) = \left | \frac{f_1(U_1(\{R_i\}_{i \in \mathcal{V}}))}{f_1(U_1(\{\tilde R_i^{-1}\}_{i \in \mathcal{V}}))}   - 1\right |.
\end{equation}
For the 
gap-function to be well-defined, we also assume that the graph $\mathcal{G}$ is connected and that the $R_{ij}$-matrices in problem $(P_1)$ are not transitively consistent. We know
according to Proposition~\ref{prop:r1:3} and Proposition~\ref{prop:r1:4} that   
\begin{align*}
& \text{gap}(\{\tilde R_i^{-1}(k)\}_{i \in \mathcal{V}}, \{\bar R_i^{-1}\}_{i \in \mathcal{V}}) \rightarrow 0, \\
& \text{gap}(\{\tilde Q_i^{-1}(k)\}_{i \in \mathcal{V}}, \{\bar R_i^{-1}\}_{i \in \mathcal{V}}) \rightarrow 0, \\
& \text{gap}(\{ Q_i(k)\}_{i \in \mathcal{V}}, \{(\text{Pr}_{O(d)}(\bar R_i))^T\}_{i \in \mathcal{V}}) \rightarrow 0, 
\end{align*}
as $k \rightarrow \infty$,
where the $\bar{R}_i$-matrices are the matrix blocks in the optimal solution $\bar{X}$ to problem $(P_2)$. By plotting the logarithm of the gap as a function of $k$, Fig.~\ref{fig:1} shows that this 
convergence holds in numerical simulations. We will now explain this figure in more detail.

The left two figures show the evolution 
of 
\begin{align*}
& \text{gap}(\{ R_i(k)\}_{i \in \mathcal{V}}, \{(\text{Pr}_{O(d)}(\bar R_i))^T\}_{i \in \mathcal{V}}) \text{ (blue line) and} \\
& \text{gap}(\{ Q_i(k)\}_{i \in \mathcal{V}}, \{(\text{Pr}_{O(d)}(\bar R_i))^T\}_{i \in \mathcal{V}}) \text{ (black line)} \hfill
\end{align*}
 for two 
different parameter settings. 
The right two figures show the evolution of
\begin{align*}
& \text{gap}(\{\tilde R_i^{-1}(k)\}_{i \in \mathcal{V}}, \{\bar R_i^{-1}\}_{i \in \mathcal{V}}) \text{ (blue line) and} \\
& \text{gap}(\{\tilde Q_i^{-1}(k)\}_{i \in \mathcal{V}}, \{\bar R_i^{-1}\}_{i \in \mathcal{V}}) \text{ (black line)} \quad \quad \quad \quad \quad 
\end{align*}
for two 
different parameter settings. The means over 100 simulations are shown. 
For each simulation $n = 10$ and $d = 5$;
the $R_{ij}$-matrices are created in a two step procedure. First $R_i$-matrices are drawn from the uniform distribution on $O(5)$. Then element-wise Gaussian noise is added with standard deviation $0.2$ to the $R_i^TR_j$-matrices. Those noisy matrices are then projected back onto $O(d)$ to obtain the $R_{ij}$-matrices. 
All the $a_{ij}$ are either equal to $1$ or $0$; the graph $\mathcal{G}$ is chosen randomly from the set of connected graphs with a certain graph density (percentage of edges). In the top figures the graph density is chosen to $0.9$ in each simulation, whereas in the 
bottom figures the graph density is chosen to $0.6$ in each simulation. The step sizes $\epsilon_1$ and $\epsilon_2$ are 
both chosen to $1/(2n)$, which is an upper bound for $\frac{2}{\|P\|_2}$.

We see that the convergence is faster when the connectivity in the 
graph, i.e., graph density, is stronger. We also see that the performance of the $R_i(k)$-matrices is better initially than that of the $Q_i(k)$-matrices,
but after sufficiently many iterations the performance of the $Q_i(k)$-matrices is several orders of magnitude better (top left and bottom left figures). This out-performance is to be expected since the gap is proven to converge to zero for the $Q_i(k)$-matrices (see Proposition~\ref{prop:r1:4}) but not for the $R_i(k)$-matrices. 

\section{Directed graphs}\label{sec:asymmetric}
In this section we introduce a distributed algorithm for the synchronization problem over connected graphs that are directed and possibly asymmetric. Such graphs are not necessarily strongly connected. The asymmetric communication scenario is present when some agents are not able to transmit information due to for example malfunctioning equipment. Asymmetric communication could also be deliberately imposed. For example if one wants a hierarchical communication structure.

\subsection{A problem formulation for directed graphs}
We begin this section by introducing the following problem
\begin{equation*}
(P_3)  
\begin{cases}
\begin{aligned}
& \underset{R}{\text{minimize}}
& & f_2(R) = \\
& & & \sum_{i = 1}^n\frac{1}{2}\| \sum_{j \in \mathcal{N}_i}{a_{ij}}(R_{ij}R_j^T - R^{T}_i) \|_F^2, \\
& \text{subject to}
& & R  = [R_1, R_2, \ldots, R_n],~\\
& & & R_i  \in O(d) \text{ for all } i.  \quad \quad \quad  \quad \quad \quad \quad  \quad \quad \quad 
\end{aligned}
\end{cases}
\end{equation*}
It is obviously not the same problem as $(P_1)$. However,
it has a structure that allows for the development of a distributed 
algorithm in situations where the graph could be asymmetric. The idea is that the optimal solution 
to this problem should be close to that of problem $(P_1)$. 
When we design an algorithm in this section the objective of the algorithm will be to calculate a good feasible solution to $(P_3)$.

We 
will now explain some connections between the two problems $(P_1)$ and $(P_3)$.
The cost function $f_2$ can be written as 
\begin{align}
f_2(R) & = \frac{1}{2}\text{tr}(U_1^TL_{\text{dir}}^TL_{\text{dir}}U_1),
\end{align} 
where 
\begin{align}
L_{\text{dir}} & = \text{diag}(A{1}_n ) \otimes I_d - W . 
\end{align} 
It holds that 
\begin{align}
L_{\text{dir}} & = L_{\text{undir}} - \text{diag}(\bar W^T \bar W) + W^T,
\end{align}
where diag in the second term on  the right-hand side is defined in same way as for 
the second term in (13).

Let $\bar{L}_{\text{dir}}$ be equivalent to the matrix ${L}_{\text{dir}}$ defined for the graph 
$\bar{\mathcal{G}}$ and the matrices $\{\bar{{R}}_{ij}\}_{(i,j) \in \bar{\mathcal{E}}}$, where $\bar{R}_{ij} = R_{ji}^T$.
It holds that
\begin{equation}
L_{\text{undir}} = {{L}_{\text{dir}} + \bar{L}_{\text{dir}}}.
\end{equation}
If the graph $\mathcal{G}$ is symmetric and $A = A^T$, it holds that
\begin{equation}
L_{\text{undir}} = {{L}_{\text{dir}} + {L}_{\text{dir}}^T}.
\end{equation}
Furthermore, if it also holds that $R_{ij}^T = R_{ji}$ for all $i$, $j$, then
\begin{equation}
L_{\text{undir}} = 2{L}_{\text{dir}}.\quad \quad \: 
\end{equation}

\subsection{The algorithm}

\rule{\columnwidth}{1pt}
\textbf{Algorithm 2} Distributed method for graphs that could be asymmetric
\hrule
\textbf{Inputs:} a directed graph $\mathcal{G} = (\mathcal{E}, \mathcal{V})$, a weight matrix $A = [a_{ij}]$, and a collection $\{R_{ij}\}_{(i,j) \in \mathcal{E}}$ of matrices in $O(d)$.\\ \\
\textbf{Outputs:}  $\tilde{R}_i(k)$, ${R}_i(k)$ for $i \in \mathcal{V}$ and $k \geq 1$.\\ \\
\textbf{Initialization:} let $\tilde{R}_i(0) = I_d$ for all $i$ and $\epsilon_3 > 0$.\\

\textbf{Iteration $k \geq 1$:} \\
for all $i$ and all $k \geq 1$ we define $\tilde{R}_i(k)$ via the following equations.
\begin{align*}
\tilde{R}_i(k) & = \tilde{R}_i(k-1) \\
& \quad + \epsilon_3\sum_{j \in \mathcal{N}_i}a_{ij}(R_{ij}\tilde{R}_j(k-1) - \tilde{R}_i(k-1)), \\
R_i^T(k) & = \text{Pr}_{\mathit{O}(d)}(\tilde{R}_i(k)).
\end{align*} 
\rule{\columnwidth}{1pt}

Let $\tilde R(k)= [\tilde R_1^T(k),\tilde R_2^T(k), \ldots, \tilde R_n^T(k)]$ for all $k$, the update for $\tilde{R}(k)$ is 
written as 
\begin{align}
\label{eq:r1:5}
\tilde{R}(k+1) & = (I - \epsilon_3 L_{\text{dir}})\tilde{R}(k).
\end{align}

\subsection{Analysis}\label{sec:zMatrix} 
The matrix $L_{\text{dir}}$ is defined for directed graphs (and asymmetric graphs). It is related to the graph connection Laplacian matrix~\cite{bandeira2013cheeger,wang2013exact,singer2012vector}, which however is only defined for undirected graphs.

\subsubsection{Some properties of $L_{\text{dir}}$}
The matrix $L_{\text{dir}}$ could be seen as a way of extending  $L_{\text{undir}}$ 
to the case of asymmetric graphs, but the matrix is also defined for symmetric graphs. We can show a strong relation between $L_{\text{undir}}$ and synchronization/transitive consistency of invertible matrices, see Proposition 8. The equivalent relation does not hold when we consider $L_{\text{dir}}$. 
The question is what we can say about $L_{\text{dir}}$. The propositions 15, 18 and 19 provide some answers to this question.

Proposition~\ref{lem:2} formulates a nullspace condition for $L_{\text{dir}}$ that must hold for transitively consistent matrices. This is a necessary condition, similar to that of Proposition~\ref{prop:2}, but under the stronger assumption that the graph is QSC. However, as stated in Proposition \ref{lem:4} further below, this condition is not a sufficient condition for synchronization of invertible matrices. Proposition 18 shows that there is an infinity of non-synchronized matrix collections close to a synchronized collection $\{R_{ij}\}_{(i,j) \in \mathcal{E}}$, all satisfying the nullspace condition. Proposition 19 states that $(P_3)$ has an optimal value equal to $0$ if and only if the $R_{ij}$'s are synchronized. It shares similarities to Lemma~\ref{lem:5} formulated for $L_{\text{undir}}$.

Before we proceed we define the matrix $\bar{L} = \text{diag}(A1_n) - A$.

\begin{prop}\label{lem:2}
For any \textbf{(QSC)} graph $\mathcal{G} = (\mathcal{V}, \mathcal{E})$,
collection $\{{R}_{ij}\}_{(i,j) \in \mathcal{E}}$ of matrices
in $GL(d,\mathbb{R})$ that is transitively consistent for $\mathcal{G}$, and collection $\{R_i\}_{i \in \mathcal{V}}$ of matrices in $GL(d,\mathbb{R})$ 
it holds that
\begin{equation}\label{eq:farfar:3}
R_{ij} = R^{-1}_{i}R_{j} \text{ for all } (i,j) \in \mathcal{E}
\end{equation}
{\textbf{(if and)} only if}
\begin{equation}\label{eq:olle}
\text{im}(\text{diag}\left(
R_{1}, R_{2}, \ldots,  R_{n} \right)V) = \text{ker}(\bar L \otimes I_d),
\end{equation}
for any matrix $V$, where the columns thereof form a basis for $\text{ker}(L_{\text{dir}})$.
In particular, if  $\mathcal{G}$ is QSC,  \eqref{eq:olle}
can be stated as 
\begin{equation}\label{eq:polle}
\text{im}(U_1(\{R_i\}_{i \in \mathcal{V}})) = \text{ker}(L_{\text{dir}}).
\end{equation}
\end{prop}

\noindent\emph{Proof: } \\
\textbf{Only if:} Suppose it holds that  
$R_{ij} = R^{-1}_{i}R_{j} \text{ for all } (i,j) \in \mathcal{E}.$
Then
\begin{align}
\label{eq:Z_laplacian}
 L_{\text{dir}} 
=~ &\text{diag}(R_1^{-1}, R_2^{-1}, \ldots, R_n^{-1})(\bar L \otimes I)\cdot \\
\nonumber
&\quad \text{diag}(R_1, R_2, \ldots, R_n).
\end{align}
Now, $L_{\text{dir}}V  = 0 \Leftrightarrow (\bar L \otimes I) \text{diag}(R_1, R_2, \ldots, R_n)V = 0 \Leftrightarrow \text{im}(\text{diag}\left(R_{1}, R_{2}, \ldots,  R_{n} \right)V)  = \text{ker}(\bar L \otimes I_d).$

\noindent 
\textbf{If:} This part only concerns the case when the 
graph $\mathcal{G}$ is QSC. 

Since $\{{R}_{ij}\}_{(i,j) \in \mathcal{E}}$ is transitively consistent 
for $\mathcal{G}$, there is a collection $\{R^{**}_i\}_{i \in \mathcal{V}}$ of matrices in $GL(d,\mathbb{R})$ such that 
\begin{align*}
L_{\text{dir}} = & \text{diag}(R_1^{**-1}, R_2^{**-1}, \ldots, R_n^{**-1})(\bar L \otimes I)\cdot \\
& \text{diag}(R^{**}_1, R^{**}_2, \ldots, R^{**}_n).
\end{align*}
Thus, the nullspace of $L_{\text{dir}}$
is given by $\text{ker}(L_{\text{dir}}) = \text{im}(V),$
where $V = \text{diag}(R^{**-1}_1, R^{**-1}_2, \ldots, R^{**-1}_n)([1, 1, \ldots, 1]^T \otimes I_d).$

Now, suppose \eqref{eq:olle} holds. Then 
$
\text{diag}\left(
R_{1}, R_{2}, \ldots,  R_{n} \right)V = ([1, 1, \ldots, 1]^T \otimes I_d)Q,$
where $Q$ is some matrix in $GL(d,\mathbb{R})$. This means that 
$$R_iR_i^{**-1}=Q \text{ for all } i,$$
which implies that $\{R^{**}_{i} \}_{i \in \mathcal{V}}$ and $\{R_{i}\}_{i \in \mathcal{V}}$ are equal up to transformation from the left. By using Lemma~\ref{lem:1} we conclude that 
$R_{ij} = R^{-1}_{i}R_{j} \text{ for all } (i,j) \in \mathcal{E}.$
\hfill $\blacksquare$\newline

\begin{rem}\label{rem:1}
In Proposition~\ref{lem:2},
the relation 
\begin{equation}
\text{im}(\text{diag}\left(
R_{1}, R_{2}, \ldots,  R_{n} \right)V) = \text{ker}(\bar L \otimes I_d)
\end{equation}
holds if and only if
for any matrix $V_2$, where the columns thereof comprise a basis for $\text{ker}(L \otimes I_d)$, there is a matrix 
$Q$
such that 
\begin{equation}
\text{diag}\left(
R_{1}, R_{2}, \ldots,  R_{n} \right)V = V_2Q.
\end{equation}
\end{rem}

\begin{rem}\label{rem:2}
In Proposition~\ref{lem:2}, if $\mathcal{G}$ is connected but not 
QSC,  it can hold that $A^T$ is the adjacency matrix of a 
QSC graph $\mathcal{G}' = (\mathcal{V}', \mathcal{E}')$. Then 
it holds that
\begin{equation}
\text{im}(U_1(\{R^{-T}_i\}_{i \in \mathcal{V}}) = \text{ker}(L_{\text{dir}}).
\end{equation} 
\end{rem}

Proposition~\ref{lem:2} is important as it states that for transitively consistent matrices and QSC graphs, the conditions \eqref{eq:1} and \eqref{eq:olle}
 are equivalent. For invertible matrices it is not 
true that \eqref{eq:olle} implies transitive consistency. See Proposition~\ref{lem:4} below. However, if we restrict all matrices to be orthogonal, the implication is true, see Proposition~\ref{prop:101} below. The latter proposition guarantees that the problem $(P_3)$ is well posed in the sense that the objective function is equal to zero if and only if the $R_{ij}$-matrices are synchronized.  

\begin{prop}\label{lem:4}
Let $\mathcal{G} = (\mathcal{V}, \mathcal{E})$ be any QSC graph
such that the node degree of one node is at least 2 (i.e, at least one element in the vector $A(\mathcal{G})[1, 1, \ldots, 1]^T$ is greater
or equal to 2). Let
 $\{{R}^*_{ij}\}_{(i,j) \in \mathcal{E}}$ be a collection of matrices
in $GL(d, \mathbb{R})$ that is transitively consistent for $\mathcal{G}$. 
Let $\{R^*_i\}_{i \in \mathcal{V}}$ be a collection of matrices in $GL(d, \mathbb{R})$ 
for which it holds that
\begin{equation}
R^*_{ij} = R^{*-1}_{i}R^*_{j} \text{ for all } (i,j) \in \mathcal{E}.
\end{equation} 

Now, for any $\epsilon > 0$, there is a collection $\{{R}_{ij}\}_{(i,j) \in \mathcal{E}}$ of matrices in $GL(d, \mathbb{R})$ that is not transitively consistent for $\mathcal{G}$ such that 
\begin{equation}\label{eq:2}
\sum_{(i,j) \in {\mathcal{E}}}\|{R}_{ij} - R^*_{ij}\|_F \leq \epsilon,
\end{equation}
and \eqref{eq:polle} holds for $\{R_{ij}\}_{(i,j) \in \mathcal{E}}$ and a collection $\{R_i\}_{i \in \mathcal{V}}$
of matrices in $GL(d, \mathbb{R})$.
\end{prop}

\noindent\emph{Proof: } \\
Suppose the $k$-th element of the vector $A(\mathcal{G})[1, 1, \ldots, 1]^T$ is larger or equal to $2$. Then there is $l,m$ such that
$l \neq k$, $m \neq k$, $R^*_{kl}, R^*_{km} \in GL(d, \mathbb{R})$.
Let $R_{kl} = R^{*-1}_{k}(I + Q)R^{*}_{l}$ and
$R_{km} = R^{*-1}_{k}(I - Q)R_m^*$, where we choose $Q \neq 0$ such that $(I \pm Q) \in GL(d, \mathbb{R})$ and 
\begin{equation}
\|Q\|_F < \frac{\epsilon}{2}.
\end{equation}
Furthermore, let 
$R_{ij} = R^*_{ij}$ for all $(i,j) \not\in \{(k,l),(k,m)\}$.
\eqref{eq:2} is satisfied.
By construction, all the $R_{ij}$ are elements of $GL(d, \mathbb{R})$.

Let $R_i = R_i^*$ for all $i$. It holds that 
\begin{align}
& L_{\text{dir}}(\mathcal{G}, \{{R}_{ij}\}_{(i,j) \in \mathcal{E}}) \\
\nonumber
=~ &\text{diag}(R_1^{-1}, R_2^{-1}, \ldots, R_n^{-1})((\bar L \otimes I) + \bar Q) \cdot \\
\nonumber
& \text{diag}(R_1, R_2, \ldots, R_n),
\end{align}
where $\mathbb{R}^{nd \times nd} \ni \bar Q = [Q_{ij}]$, $Q_{kl} = a_{kl}Q$, $Q_{km} = -a_{km}Q$ and $Q_{ij} = 0 \in \mathbb{R}^{d \times d}$ for 
all $(i,j) \not\in \{(k,l),(k,m)\}$. It holds that $\text{ker}((\bar L \otimes I) +\bar Q)
\supset \text{ker}(\bar L \otimes I)$. Since the eigenvalues of $L_{\text{dir}}(\mathcal{G}, \{{R}_{ij}\}_{(i,j) \in \mathcal{E}})$ are continuous functions of $Q$, it holds that for $Q$ small enough, $\text{ker}((\bar L \otimes I) + \bar Q)
\subset \text{ker}(\bar L \otimes I)$. We assume without loss of generality that $Q$ has been chosen small enough so that the latter inclusion holds.
Then \eqref{eq:olle} holds for the $R_i$-matrices.
According to Proposition~\ref{lem:2}, if the $R_{ij}$-matrices are transitively consistent and $\mathcal{G}$ is QSC, \eqref{eq:olle} is a condition to guarantee \eqref{eq:1}. But \eqref{eq:1} is not fulfilled 
since $R_kR_{kl}R_l^{-1} = (I + a_{kl}Q) \neq I$. Thus, the $R_{ij}$-matrices are
not transitively consistent.
\hfill $\blacksquare$\newline

\begin{prop}\label{prop:101}
For any QSC-graph $\mathcal{G} = (\mathcal{V}, \mathcal{E})$ and 
collection $\{{R}_{ij}\}_{(i,j) \in \mathcal{E}}$ of matrices
in $O(d)$ 
it holds that  $\{{R}_{ij}\}_{(i,j) \in \mathcal{E}}$ is transitively consistent 
if and only if
there is a  collection $\{R_i\}_{i \in \mathcal{V}}$ of matrices in $O(d)$ such that 
\begin{equation}\label{eq:polle:300}
\text{im}(U_1(\{R_i\}_{i \in \mathcal{V}})) \subset \text{ker}( L_{\text{dir}}).
\end{equation}
\end{prop}

\noindent\emph{Proof: }\\
\textbf{Only if:}  Suppose that the conditions in Proposition~\ref{prop:101} on the graph and the matrices hold. Suppose $\{{R}_{ij}\}_{(i,j) \in \mathcal{E}}$ is transitively consistent. There is a collection $\{R_i\}_{i \in \mathcal{V}}$ of matrices in $O(d)$ such that \eqref{eq:farfar:3} holds for the $R_i$-matrices. According to Proposition~\ref{lem:2}, the condition \eqref{eq:polle} also holds for the $R_i$-matrices.

\textbf{If:} 
Suppose $\{{R}_{ij}\}_{(i,j) \in \mathcal{E}}$ is not transitively consistent and that there is a collection $\{R_i\}_{i \in \mathcal{V}}$ of matrices in $O(d)$ such that \eqref{eq:polle:300} holds. Let $R_i(R_{ij})R_j^T = (I + Z_{ij})$ for all $(i,j) \in \mathcal{E}$. We know that there is a least one pair $(i_1,j_1) \in \mathcal{E}$ such that $Z_{i_1j_1} \neq 0$ (otherwise  $\{{R}_{ij}\}_{(i,j) \in \mathcal{E}}$ would be transitively consistent).

Now, 
\begin{align}
& L_{\text{dir}}(\mathcal{G}, \{{R}_{ij}\}_{(i,j) \in \mathcal{E}}) \\
\nonumber
=~ &\text{diag}(R_1^{T}, R_2^{T}, \ldots, R_n^{T})((\bar L \otimes I) + \bar Z) \cdot \\
\nonumber
& \text{diag}(R_1, R_2, \ldots, R_n),
\end{align}
where $\bar{Z} \in \mathbb{R}^{nd \times nd}$ is a block matrix whose $(i,j)$-th block is equal to $a_{ij}Z_{ij} \in \mathbb{R}^{d \times d}$. 

It holds that
$L_{\text{dir}}U_1(\{R_i\}_{i \in \mathcal{V}})$ is equal to  zero if and only if $\bar Z[I, I, \ldots, I]^T$ is equal to zero. Thus,  
it holds that 
\begin{equation}\label{eq:farfar:1}
\sum_{j \in \mathcal{N}_{i_1}}a_{i_1j}Z_{i_1j} = 0.
\end{equation}
But it also holds, since the $(I + Z_{ij})$-matrices are orthogonal,  that
\begin{align}
\label{eq:farfar:2}
Z_{i_1j}^TZ_{i_1j} + Z_{i_1j} + Z_{i_1j}^T = 0, \text{ for all } j \in \mathcal{N}_{i_1}.
\end{align}
By combining \eqref{eq:farfar:1} and \eqref{eq:farfar:2} we obtain that 
$$\sum_{j \in \mathcal{N}_{i_1}}a_{i_1j}Z_{i_1j}^TZ_{i_1j} = 0,$$
which implies (since all the $a_{ij} > 0$) that $Z_{i_1j_1} = 0$. But this is a contradiction to the statement that $Z_{i_1j_1} \neq 0$.
\hfill $\blacksquare$\newline

The difference between the two previous proofs is
the following.
In the proof of Proposition~\ref{lem:4} we change $L_{\text{dir}}$ with additional matrices such that the nullspace is preserved. In the proof of Proposition~\ref{prop:101}, when we try to preserve the nullspace of $L_{\text{dir}}$ under the constraint that the changed matrices are orthogonal, we reach the conclusion that the change of the matrices must be equal to zero. This means that we have to go outside the set of orthogonal matrices when we change $L_{\text{dir}}$.

\subsubsection{Convergence analysis}\label{sec:r1:conv2}
\begin{table}[!ht]
  \begin{center}
    \begin{tabular}{| l | l |}
    \hline
    (1) & $\mathcal{G}= (\mathcal{V}, \mathcal{E})$ is QSC.\\
    \hline
    (2) & $\mathcal{G}= (\mathcal{V}, \mathcal{E})$ is strongly connected.\\
        \hline
    (3) & $R_{ij} \in O(d)$ for all $(i,j) \in \mathcal{E}$.\\
        \hline
    (4) & $\{R_{ij}\}_{(i,j) \in \mathcal{E}}$ is transitively consistent.\\
        \hline
    (5) & Let $\text{Re}(\lambda_1) \leq \text{Re}(\lambda_2) \leq \cdots \leq  \text{Re}(\lambda_n)$, \\
    & where the $\lambda_i$'s are the eigenvalues \\
    & of $\bar L = \text{diag}(A1_n) - A$. It holds that \\
    &  $|1 - \epsilon_3 \lambda_i| < 1$  for $i = 1, 2, \ldots, n-1$. \\
    \hline
    (6) & There are $\delta_1, \delta_2 > 0$ such that $\delta_1 \leq a_{ij} \leq \delta_2$ \quad \quad \quad \quad \\
    &  for all the $a_{ij}$-weights. \\
    \hline
    \end{tabular}
  \end{center}
  \caption{Conditions for convergence.}
\end{table}
 
Now we introduce the following result, in analogy to Proposition~\ref{prop:r1:2}, but defined 
for asymmetric graphs. 
\begin{prop}
Suppose convergence conditions (1) and (3-5) hold.
Then, for Algorithm 2, there is a positive integer $K$ such that $\tilde{R}_j^{-1}(k)$ is well-defined for $k\geq K$, and (for $k \geq K$) it holds that
$$(R_i^T(k)R_j(k) \rightarrow R_{ij} \text{ as } k \rightarrow \infty) \text{ for all } (i,j) \in \mathcal{E},$$
and 
$$(\tilde{R}(k)\tilde{R}_j^{-1}(k) \rightarrow R_{ij} \text{ as } k \rightarrow \infty) \text{ for all } (i,j) \in \mathcal{E}.$$
\end{prop}

\noindent\emph{Proof: } 
When conditions (3) and (4) hold it is easy to show that $L_{\text{dir}}$ is similar to $\bar{L} \otimes I_d$, where $\bar{L}$ is defined in condition (5).
Since conditions (1) and (3) hold, we can use Proposition~\ref{lem:2}. There is a collection $\{\bar R_i\}_{i \in \mathcal{V}}$ such that 
$$\text{ker}(L_{\text{dir}}) = \text{im}(U_1(\{\bar R_i\}_{i \in \mathcal{V}})).$$
The rest of the proof is equivalent to that of Proposition~\ref{prop:r1:2} and omitted.
\hfill $\blacksquare$
 
For the general case when transitive consistency does not hold, the  convergence results in this section are more 
conservative than those provided in Section~\ref{sec:r1:conv}
for Algorithm 1. The main difficulty is that the matrix $L_{\text{dir}}$ is not symmetric, and we cannot guarantee in 
general that its spectrum is real.
We have the following theoretical result.

\begin{prop}\label{lem:10}
Suppose convergence conditions (1), (3), and (6) hold.
Then 
\begin{enumerate}
\item The dynamical system 
\begin{equation*}
\dot{x}(t) = -{L_{\text{dir}}}x(t), \quad x(t) \in \mathbb{R}^{nd},
\end{equation*}
is critically stable.\\
\item $\|P\|_2$ is an upper bound for
$\|L_{\text{dir}}\|_2$ 
where $P = \text{diag}(A1_n) + A$. \\
\item If also the convergence condition (2) holds,  \\
\begin{enumerate}
\item the matrix $L_{\text{dir}}$ has no purely imaginary eigenvalues, \\
\item there is $\epsilon(n,d, \delta_1, \delta_2) > 0$, such that for $\epsilon_3 \in (0,\epsilon)$
the eigenvalues of the matrix $(I - \epsilon_3 L_{\text{dir}}(\mathcal{G}, \{{R}_{ij}\}_{(i,j) \in \mathcal{E}}))$ are contained in the 
closed unit disc and if there is an eigenvalue on the boundary 
of the disc, it is equal to $1$ and the algebraic respective
geometric multiplicities thereof are equal. 
\end{enumerate}

\end{enumerate}
\end{prop}

\noindent\emph{Proof: } \\
\emph{(1)} Let 
\begin{equation}\label{eq15}
\dot{x}(t) = -{L_{\text{dir}}}x(t), \quad x(t) \in \mathbb{R}^{nd},
\end{equation}
where $x(0)$ is the initial state.
We can write $x(t)$ as $x(t) = [x_1^T(t), x_2^T(t), \ldots, x_n^T(t)]^T$, where $x_i(t) \in \mathbb{R}^{d}$ for all $i$. 
Define the function 
$V(x) = \max_i(x_i^Tx_i).$
If there is some eigenvalue  of ${L_{\text{dir}}}$
with negative real part or if there is a Jordan 
block of dimension larger than one corresponding to an eigenvalue pair on the imaginary axis, there is $x_0$ such that for the state $x(t)$ with initial
state $x_0$, $V(x(t)) \rightarrow \infty$ as $t \rightarrow \infty$. We want to show that this is not possible. Let us first define the set
$\mathcal{I}_{\max}(t) = \{i:V(x(t)) = x_i^T(t)x_i(t)\}.$
Now,
\begin{align}
\label{eq30}
D^+(V(x(t))) & = \max_{i \in \mathcal{I}_{\max}(t)}\frac{d}{dt}x_i^T(t)x_i(t) \\
\nonumber
& =  \max_{i \in \mathcal{I}_{\max}(t)}x_i^T(t)\left(\sum_{j \in \mathcal{N}_i}({R}_{ij}x_j(t) - x_i(t))\right ) \\
\nonumber
& \leq 0,
\end{align}
where $D^+$ is the upper Dini-derivative. A proof 
of the first equality in~\eqref{eq30} can be found in~\cite{thunberg2014consensus}
using the results in \cite{yoshizawa1966stability} and \cite{clarke1975generalized}.   
Now we can use the 
Comparison Lemma~\cite{khalil2002nonlinear} to show that $V(x(t))$ is decreasing independently 
of the choice of $x_0$. The last inequality in~\eqref{eq30}
is a consequence of the fact that the ${R}_{ij}$-matrices
are orthogonal.

\emph{(2)} 
Let $v = [v_1^T, v_2^T, \ldots, v_n^T]^T \in \mathbb{S}^{nd}$, where
each $v_i \in \mathbb{R}^{d}$. By using the structure of the $f_2$-function,
one can show that 
$v^TL_{\text{dir}}^TL_{\text{dir}}v  = \sum_{i = 1}^n\| \sum_{j \in \mathcal{N}_i}{a_{ij}}(R_{ij}v_j - v_i) \|_F^2.$
Now, 
\begin{align}
& v^TL_{\text{dir}}^TL_{\text{dir}}v  \leq \sum_{i = 1}^n\| \sum_{j \in \mathcal{N}_i}{a_{ij}}(\|v_j\|_2 + \|v_i\|_2) \|_2^2 \\
\nonumber
\leq & \|P\|^2_2.
\end{align}

\emph{(3)}  Now we show, by a contradiction argument, that there are no purely imaginary eigenvalues. Suppose there are purely imaginary eigenvalues. Then there must be a nontrivial periodic solution
$\bar{x}(t) = [\bar{x}_1^T(t), \bar{x}_1^T(t), \ldots, \bar{x}_n^T(t)]^T$ to  \eqref{eq15}, i.e., $\bar{x}(t)$ is periodic and $\bar{x}(t) \neq \bar{x}(t + \Delta)$ for all $t$ and some 
$\Delta > 0$.

 It can be shown that in this particular situation $D^+(V(\bar{x}(t))) = 0$ for all $t$ and it can also be shown that a necessary condition for this to hold is that $\bar{x}_i(t) = \bar{x}_j(t)$ for all $i,j$ and $t$, and furthermore, this can only happen if $R_{ij} = I$ for all $(i,j)$. 
 
In the following three paragraphs we prove that the claims in the preceding paragraph are true. If $D^+(V(\bar{x}(t_1))) \neq 0$, it holds that (see the proof of part \emph{(1)}) $D^+(V(\bar{x}(t_1))) < 0$. But then 
$V(\bar{x}(t)) < V(\bar{x}(t_1))$ for all $t > t_1$ (see the proof of part \emph{(1)}), which implies
that there is no $t_2 > t_1$ such that $\bar{x}(t_2) = \bar{x}(t_2)$. Thus the solution is not periodic, which is a contradiction. Hence, we can conclude that $D^+(V(\bar{x}(t))) = 0$ for all $t$. 

Now, let us assume that there is $i_1,j_1 \in \mathcal{V}$ such that $\bar{x}_{i_1}(t_1) \neq \bar{x}_{j_1}(t_1)$ at some time $t_1$. Let
$V_{\max}(t_1,t) = \{k:\|x_k(t)\| = \max_{l \in \mathcal{V}}\|x_l(t_1)\|\}.$ We can without loss of generality assume that
$i_1 \in V_{\max}(t_1,t_1)$,
$j_1 \not\in V_{\max}(t_1,t_1)$, and that $j_1 \in \mathcal{N}_{i_1}$.  The reason for why we can make this assumption about $i_1$ and $j_1$, is that the graph $\mathcal{G}$ is strongly connected. 
Now it holds that
\begin{align}
\nonumber
& \frac{d}{dt}x_{i_1}^T(t)x_{i_1}(t)|_{t_1} \\
\label{eq:nisse:200}
=&x_{i_1}^T(t_1)\sum_{j \in \mathcal{N}_{i_1}}({R}_{{i_1}j}x_j(t_1) - x_{i_1}(t_1)) < 0.
\end{align}
The strict inequality holds because $x_{i_1}^T(t_1)({R}_{{i_1}j_1}x_{j_1}(t_1) - x_{i_1}(t_1)) < 0$.
We also know that $\bar{x}(t)$ is continuous and
continuously differentiable. There must be $t_2 > t_1$
such that $V_{\max}(t_1,t_2) \subset V_{\max}(t_1,t_1) - \{i_1\}.$ 
Now, assume that $V_{\max}(t_1,t_2)$ is nonempty. Then we
can perform the analogous procedure again, where we pick $i_2 \in V_{\max}(t_1,t_2)$,
$j_2 \not\in V_{\max}(t_1,t_2)$, where $j_2 \in \mathcal{N}_{i_2}$, and show that that there is a $t_3 > t_2$ such that $V_{\max}(t_1,t_3) \subset V_{\max}(t_1,t_2) - \{i_2, i_1\}.$ By repeating this procedure, one can show by induction that there must be a $t_f > t_1$ such that $V_{\max}(t_1,t_f) = \emptyset$. But  $V(\bar{x}(t_f)) < V(\bar{x}(t_1))$ contradicts the fact that $\bar{x}(t_f) = \bar{x}(t_1)$.

To show that $\bar{x}_i(t) = \bar{x}_j(t)$ for all $i,j$ and $t$ can only happen when $R_{ij} = I$ for all $(i,j)$, we note the following. Suppose it holds that $\bar{x}_i(t) = \bar{x}_j(t)$ for all $i,j$ and $t$ and in particular at the time $t_1$. Suppose there is $i$ and an $j \in \mathcal{N}_i$ 
such that $R_{ij} \neq I_d$. But then one can show that the time derivative of $\|\bar x_i(t_1)\|^2$ is negative at the time $t_1$. To see this, consider the structure of the right-hand side of \eqref{eq:nisse:200}; it holds that $\bar x_{i}^TR_{ij}\bar x_j - \bar x_i^T\bar x_i < 0$.
 Thus there is a time $t_2 > t_1$ where $\|\bar x_i(t_2)\|^2 < \|\bar x_i(t_1)\|^2$. On the other hand, since  $D^+(V(\bar{x}(t))) = 0$ for all $t$, there must be a $k \in \mathcal{V}$ such that $\|\bar x_k(t_2)\|^2 = \|\bar x_i(t_1)\|^2$. But then $\|\bar x_k(t_2)\|^2 \neq \|\bar x_i(t_2)\|^2$, which contradicts the fact that $\bar x_i(t_2) = \bar x_k(t_2)$ holds.

Now, if the $R_{ij} \neq I$, the 
necessary condition is not fulfilled, hence we have a contradiction. In the case when the $R_{ij} = I$, it holds that $L_{\text{dir}}(\mathcal{G}) = \bar L \otimes I_d$ and the latter matrix 
does not have purely imaginary eigenvalues.

The $R_{ij}$-matrices and the $a_{ij}$-weights are contained in a compact set parameterized by $(n,d, \delta_1, \delta_2)$, see convergence condtion (6). The eigenvalues of $L_{\text{dir}}$ are continuous over this set. 
As setforth above, there are no purely imaginary eigenvalues. Hence, any complex eigenvalue must have an argument in the open set $(\pi/2, 3\pi/2)$. 
Supremums and infimums of continuous functions over compact sets amount to maximums and minimums, respectively, see Weierstrass' Extreme Value Theorem.
This means that the arguments of the complex eigenvalues are contained in a set $[\pi/2+ \tilde \delta_1, 3\pi/2 - \tilde \delta_2]$ for some $\tilde{\delta}_1, \tilde{\delta}_2  > 0$. Furthermore,  with analogous reasoning, the magnitudes of the eigenvalues are bounded from above by some $\tilde{\delta}_3$ (the magnitudes are continuous functions). 
Now we can construct $\epsilon$. 
Let $\alpha_1$ and $\alpha_2$ be the unique positive solutions to
\begin{align*}
\|[\alpha \cos(\pi/2 + \tilde{\delta}_1) - 1, \alpha \sin(\pi/2 + \tilde{\delta}_1)]^T\|  & = 1 \text{ and } \\
\|[\alpha \cos(3\pi/2- \tilde \delta_2) - 1, \alpha \sin(3\pi/2- \tilde \delta_2)]^T\|  & = 1,
\end{align*}
respectively.
We choose $\epsilon = \min\{\alpha_1, \alpha_2\}/\tilde \delta_3$.

Now, any eigenvalue that is equal to $1$ for $(I - \epsilon_3L_{\text{dir}})$ is corresponding to an eigenvalue that is equal to $0$ for $L_{\text{dir}}$. The eigenvalues have been shifted by $1$ between the two. The eigenvectors are the same.  
If the claim in \emph{(3)} about the multiplicity would be false, then $L_{\text{dir}}$ would be unstable, which is a contradiction, see \emph{(1)}.
\hfill $\blacksquare$\newline

Suppose the graph
$\mathcal{G}$ is QSC and the collection of matrices $\{R_{ij}\}_{(i,j) \in \mathcal{E}}$ is transitively consistent. From Proposition~\ref{lem:10} we can draw the following conclusion about Algorithm 2.  If $\epsilon_3$ is chosen small enough, the system defined in \eqref{eq:r1:5} is stable. However, we do not provide any explicit bounds for $\epsilon_3$. An interesting research direction to pursue would be to investigate if the results in~\cite{barrau2015intrinsic} can be used to derive such bounds.

\subsection{Numerical simulations}\label{sec:simulations}
\begin{figure}[!th] 
\centering
\includegraphics[scale=0.23]{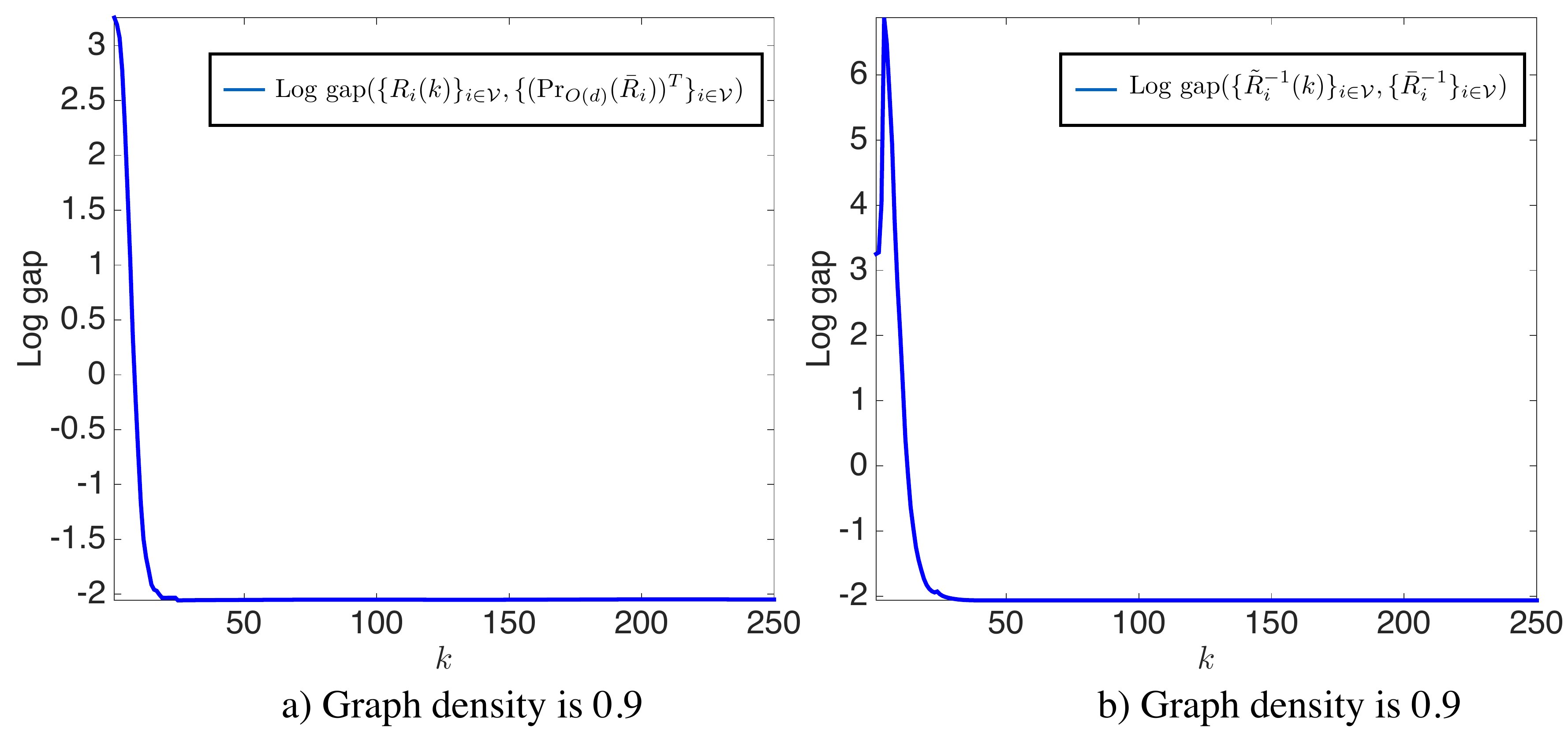}
\centering
\includegraphics[scale=0.23]{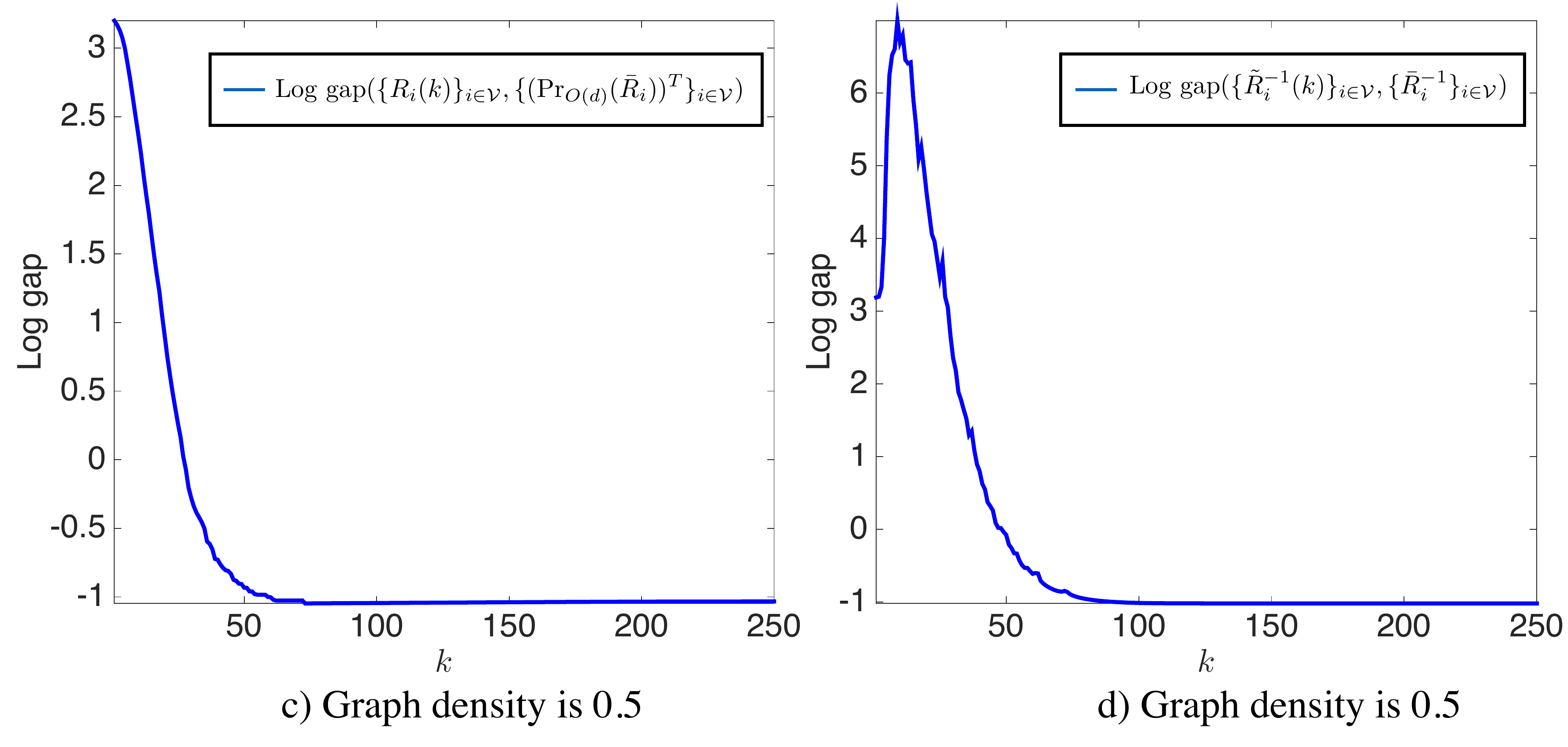}
\caption{Convergence plots for Algorithm 2.}
 \label{fig:2}
\end{figure}

In Figure~\ref{fig:2} the evolution of 
\begin{align*}
& \text{gap}(\{ R_i(k)\}_{i \in \mathcal{V}}, \{(\text{Pr}_{O(d)}(\bar R_i))^T\}_{i \in \mathcal{V}}) \text{ (left two figures),} \\
& \text{gap}(\{\tilde R_i^{-1}(k)\}_{i \in \mathcal{V}}, \{\bar R_i^{-1}\}_{i \in \mathcal{V}}) \text{ (right two figures)}
\end{align*}
is shown. The definition of the gap function and the problem setup are the same as in Section~\ref{sec:num:1}, apart from the construction of the graphs.
The graphs are now QSC and in general asymmetric. Moreover, $\epsilon_3 = 1/(2n)$. For the top two figures the graph density is $0.9$ and for the bottom two figures the graph density is
$0.5$.

There are two important things to note here. None of the results
in the previous section guarantee that the algorithm converges when the graph is QSC and not necessarily 
strongly connected, nor do they guarantee that the algorithm converges when the step size is $1/(2n)$. Yet this was the case in all simulations we conducted. We let the graph density vary 
between $0.3$ and $0.9$, the dimension vary between $3$ and $20$, 
we tried different noise levels and $a_{ij}$-weights, and convergence to at least $-0.5$ for the gap (in logarithmic scale) was observed for the means.
 
\section{Conclusion}\label{sec:conclusion}
In this work we proposed two distributed methods to synchronize orthogonal matrices of  arbitrary dimension over  graphs. One was presented for symmetric graphs and one was presented for graphs that could also be asymmetric. For both methods we provided a theoretical analysis including various convergence conditions. In the case of symmetric graphs, we proved convergence to the solution of a spectral relaxation. Numerical simulations verify this result. In the case of directed and asymmetric graphs, we proved convergence of the algorithm for small step sizes. This convergence was illustrated in numerical simulations.

\bibliographystyle{unsrt}       
\bibliography{the_paper.bbl}
\end{document}